\documentclass[11pt]{article}
\usepackage{latexsym, amssymb, amsfonts, amsmath, amsthm, graphics, graphicx, subfigure, bbm, stmaryrd}
\usepackage{epstopdf}
\newcommand{\field}[1]{\mathbb{#1}}
\newcommand{\R}{\field{R}}
\newcommand{\E}{\field{E}}
\newcommand{\C}{\field{C}}
\newcommand{\Z}{\field{Z}}

\newcommand{\Tr}{\mathrm{Tr}}
\newcommand{\tr}{\mathrm{tr}}
\newcommand{\cov}{\mathrm{cov}}
\newcommand{\var}{\mathrm{var}}
\newtheorem{theorem}{Theorem}

\newtheorem{proposition}{Proposition}

\newtheorem{lemma}{Lemma}
\newtheorem{corollary}{Corollary}
\newtheorem{remark}{Remark}

\newcommand{\met}[2][ccccccccccccccccccccccccccccccccc]{\left[ \begin{array}{#1} #2 \\ \end{array}\right]}

\numberwithin{equation}{section}

\title{\bf{General Tridiagonal Random Matrix Models, \\ Limiting Distributions and Fluctuations }}
\author{ Ionel Popescu \\  \\ School of Mathematics, Georgia Institute of Technology, 686 Cherry Street, Atlanta GA, 30332, USA \\ and   \\ IMAR, 21 Calea Grivitei Street, 010702-Bucharest, Sector 1, Romania \\ email: \texttt{ipopescu$@$math.gatech.edu}}
\date{}

\begin{document}
 \maketitle

\abstract{In this paper we discuss general tridiagonal matrix models which are natural extensions of the ones given in \cite{DE} and \cite{DE2}.  We prove here the convergence of the distribution of the eigenvalues and compute the limiting distributions in some particular cases.  We also discuss the limit of fluctuations, which, in a general context, turn out to be Gaussian.  For the case  of several random matrices, we prove the convergence of the joint moments and the convergence of the fluctuations to a Gaussian family.  

The methods involved are based on an elementary result on sequences of real numbers and a judicious counting of levels of paths. }

\section{Introduction}

Tridiagonalization is a standard procedure in numerical analysis.  The advantage of  tridiagonalization is  that the eigenvalues do not change under this procedure on one hand and on the other hand the tridiagonal matrix is easier to study, both numerically and theoretically.

The well known GOE, GUE and GSE random matrix models (see \cite{M} for a standard reference), have the eigenvalue distribution given by the density
\begin{equation}\label{e00}
\frac{1}{Z_{n,\beta}}\prod_{1\le i<j\le n}|x_{i}-x_{j}|^{\beta}e^{-\beta\sum x_{i}^{2}/2},
\end{equation}
for $\beta=1,2,4$ and $Z_{n,\beta}$ is the corresponding normalization constant.  

For $\beta=2$, tridiagonalizing the GUE ensembles, in \cite{DE} and \cite{DE2}, the authors arrive at
\begin{equation}\label{e0}
\frac{1}{\sqrt{\beta}}\met{N(0,2)& \chi_{(n-1)\beta }&0&0&0&\dots & 0 & 0 &0 \\  \chi_{(n-1)\beta} & N(0,2) & \chi_{(n-2)\beta} &0 &0 &\dots &0&0 &0\\ 0& \chi_{(n-2)\beta} & N(0,2) & \chi_{(n-3)\beta} & 0 & \dots &0&0&0 \\ \hdotsfor{9}\\ \hdotsfor{9} \\ 0&0 & 0 & 0 & 0& \dots & \chi_{2\beta} & N(0,2) & \chi_{\beta}  \\ 0 & 0 & 0& 0&0& \dots &0 & \chi_{\beta} & N(0,2) }
\end{equation}
where all entries are independent and $\chi_{r}$ is the $\chi$ distribution with $r$ degrees of freedom.  Since the tridiagonalization does not change the eigenvalue distribution, it follows that for this model the eigenvalues have the distribution given by \eqref{e00}.  Moreover it turns out that for any arbitrary $\beta>0$, the eigenvalue distribution of the model \eqref{e0} is given by \eqref{e00}.  

Obviously the models \eqref{e0} are less complex and consequently one should be able to take advantage of this,  particularly in the case of computations of expectations of traces of powers.   In  \cite{DE} and \cite{DE2} the limit distribution and the fluctuations are studied.  However some of the arguments used there rely on the particular form of the model and it's not clear weather these particular properties are really needed for the convergence and fluctuations.  

Another model which is discussed in the literature is the Wigner ensemble which appeared for the first time in \cite{W1} and \cite{W2}.   These are symmetric random matrices with upper diagonal entries independent of one another with mean zero and the same variance.  For these ensembles,  Wigner himself proved a form of convergence of the distribution of eigenvalues to the semicircle law.  The main method available here to study the limiting eigenvalue distribution and fluctuations is so called moment method which consists in expanding the traces of powers and counting the contributing terms.  There are various sources using this method, among the so many we mention for instance \cite{SS}  and the survey paper  \cite{B} for various combinatorial but also analytic approaches.  For the problem of fluctuations from the limiting distribution, a very general form can be found in \cite{AZ}.  Another use of the moment problem is in \cite{SO2} for universality at the edge of the spectrum.  

In the context of tridiagonal models we would like to introduce and discuss the analog of the Wigner ensembles and prove the convergence of the distribution of the eigenvalues and the fluctuations using the method of moments.  We show a nice and clean combinatorial way of doing this.  

At first, these models may seem to be an extension in form only.  There are many reason we want to study these.  The first one is that these seem to be the natural analog of the Wigner ensembles for the tridiagonal ensembles.  It turns out that these ensembles obey nice properties as convergence of the empirical distribution of the eigenvalues and the fluctuations converging to a Gaussian family.  Thus these can be seen as another universality property.  The second reason, the main one is that tridiagonalization of a Wiegner ensemble outputs a random tridiagonal matrix.  We are still far from understanding these tridiagonal matrices due to the fact that the entries of the resulting matrix are no longer independent.  What we try here is to study models in which the entries are independent with the hope that these will shed light on the more intricate case with dependencies.  The third reason is connected to the following problem. Take a band matrix of width growing with the size of the matrix.  These models have been studied in the literature in some situations, but there are various cases where not much is known.  Such case is the one in which the band width is the square root of the matrix size.  To the knowledge of the author it appears that the convergence of the empirical distribution is not known.  The tridiagonalization of such an ensemble produces a matrix whose entries have strong dependencies but we believe that studying these will bring to light some interesting phenomena.

Our main matrix model in this paper is given by 
\begin{equation}\label{def}
A_{n}=\met{d_{n}& b_{n-1}&0&0&0&\dots & 0 & 0 &0 \\ b_{n-1} & d_{n-1} & b_{n-2} & 0 &0 &\dots &0&0 &0\\ 0& b_{n-2} & d_{n-2} & b_{n-3} & 0 & \dots &0&0&0 \\ \hdotsfor{9}\\ \hdotsfor{9} \\ 0&0 & 0 & 0 & 0& \dots & b_{2} & d_{2} & b_{1}  \\ 0 & 0 & 0& 0&0& \dots &0 & b_{1} & d_{1} }
\end{equation}
where the entries are independent random variables.  In particular if $\{d_{n}\}_{n=1}^{n}$ is a sequence of iid normal random variables and $b_{n}=\chi_{n\beta}/\sqrt{\beta}$, then we get \eqref{e0}.  One of the main properties used in \cite{DE} and \cite{DE2} to study the limiting eigenvalue distribution and the fluctuation is the simple fact that $\chi_{r}-\sqrt{r}$ converges in distribution to $N(0,1/2)$.  Rephrased, it implies that in distribution sense
\begin{equation}\label{e:1}
\lim_{n\to\infty} b_{n}/\sqrt{n}=1.   
\end{equation}
This together with the fact that $d_{n}$ are iid with finite moments, turn out to be sufficient for proving the convergence of the eigenvalues to the semicircle law for the rescaled matrix $X_{n}=\frac{1}{\sqrt{n}}A_{n}$. 

In what follows, for any matrix $Y=(y_{i,j})_{i,j=1}^{n}$, we use $\Tr_{n}(Y)=\sum_{i=1}^{n}y_{i,i}$ for the full trace and $\tr_{n}(Y)=\frac{1}{n}\sum_{i=1}^{n}y_{i,i}$ for the reduced trace. 

To outline the idea of this paper in one instance, namely the convergence in moments of the eigenvalue distribution,  let's take the trace of the fourth moment of $X_{n}$, which is 
\[
\tr_{n}(X_{n}^{4})=\frac{1}{n^{3}}\sum_{1\le i_{1},i_{2},i_{3},i_{4}\le n}a_{i_{1},i_{2}}a_{i_{2},i_{3}}a_{i_{3},i_{4}}a_{i_{4},i_{1}}.
\]
We want to show that this converges.  Here $a_{i,j}$ are the entries of the matrix $A_{n}$.  Now since the matrix $A_{n}$ is tridiagonal, these terms are zero for       
 $|i_{u}-i_{u+1}|\ge2$ for $1\le u\le 4$ with  $i_{5}=i_{1}$.  Hence the only nonzero contribution is given by the sequences $(i_{1},i_{2},i_{3},i_{4})$ with $|i_{u}-i_{u+1}|\le 1$.  Let's  call these sequences admissible.   Now we rewrite 
\begin{equation}\label{eq1}
\tr_{n}(X_{n}^{4})=\frac{1}{n^{3}}\sum_{p=1}^{n}\sum_{\substack{1\le i_{1},i_{2},i_{3},i_{4}\le n \:\text{admissible} \\ \max(i_{1},i_{2},i_{3},i_{4})=p } }a_{i_{1},i_{2}}a_{i_{2},i_{3}}a_{i_{3},i_{4}}a_{i_{4},i_{1}},
\end{equation}
Since the  indices $i_{1},i_{2},i_{3},i_{4}$ are in within finite distance from one another, for $p$ larger than $3$, the sum 
\[
S_{p}=\sum_{\substack{1\le  i_{1},i_{2},i_{3},i_{4}\le n \:\text{admissible} \\ \max(i_{1},i_{2},i_{3},i_{4})=p } }a_{i_{1},i_{2}}a_{i_{2},i_{3}}a_{i_{3},i_{4}}a_{i_{4},i_{1}}=\sum_{\substack{p-2\le i_{1},i_{2},i_{3},i_{4}\le p \:\text{admissible} \\ \max(i_{1},i_{2},i_{3},i_{4})=p } }a_{i_{1},i_{2}}a_{i_{2},i_{3}}a_{i_{3},i_{4}}a_{i_{4},i_{1}}
\]
depends only on $p$ and not on $n$ and 
\[
\E[\tr_{n}(X_{n}^{4})]=\frac{1}{n^{3}}\sum_{p=1}^{n}\E[S_{p}].
\]
In the limit, for large $n$, one can ignore from the above sum, the terms $\E[S_{1}]$, $\E[S_{2}]$ and $\E[S_{3}]$ or any finite number of them.  The key to our computations is the following simple result on sequences:
\[
\lim_{p\to\infty}\frac{x_{p}}{p^{2}}=M \Longrightarrow \lim_{n\to\infty}\frac{1}{n^{3}}\sum_{p=3}^{n}x_{p}=\frac{M}{3}.
\]
Applying this to $x_{p}=\E[S_{p}]$, one reduces the computation of the limit of $\E[\tr_{n}(X_{n}^{k})]$ to 
\[
\lim_{p\to\infty}\frac{1}{p^{2}} \E[S_{p}].
\]
Next, we notice that $i_{1}=j_{1}+p$, $i_{2}=j_{2}+p$, $i_{3}=j_{3}+p$ and $i_{4}=j_{4}+p$, with $|j_{u}-j_{u+1}|\le 1$, $1\le u\le 4$, $j_{5}=j_{1}$ and rewrite 
\[
S_{p}=\sum_{\substack{ j_{1},j_{2},j_{3},j_{4}\le 0 \:\text{admissible} \\ \max(j_{1},j_{2},j_{3},j_{4})=0 } }a_{j_{1}+p,j_{2}+p}a_{j_{2}+p,j_{3}+p}a_{j_{3}+p,j_{4}+p}a_{j_{4}+p,j_{1}+p}
\]
and the limit of $\E[S_{p}]/p^{2}$ reduces to the ones of the form
\[\tag{*}
\lim_{p\to\infty}\frac{1}{p^{2}}\E[a_{j_{1}+p,j_{2}+p}a_{j_{2}+p,j_{3}+p}a_{j_{3}+p,j_{4}+p}a_{j_{4}+p,j_{1}+p}].
\]
If $j_{u}=j_{u+1}$, then $a_{j_{u}+p,j_{u+1}+p}=d_{j_{u}+p}$, while for the case $|j_{u}-j_{u+1}|=1$, we have $a_{j_{u}+p,j_{u+1}+p}=b_{j_{u}+p}$, hence if at least one of the situations $j_{u}=j_{u+1}$ occur, \eqref{e:1} forces the limit in (*) to be $0$.  This means that the only contributing terms are those with all consecutive $j$'s different.  This happens if and only if $(j_{1},j_{2},j_{3},j_{4},j_{1})$ is one of the following 
\begin{align}\label{gs}
\notag\gamma_{1}&=(0,-1,0,-1,0), \: & \gamma_{2}&=(0,-1,-2,-1,0),\: &\gamma_{3}&=(-1,0,-1,0,-1),\\ \gamma_{4}&=(-1,-2,-1,0,-1),\: &\gamma_{5}&=(-1,0,-1,-2,-1),\:&\gamma_{6}&=(-2,-1,0,-1,-2).
\end{align}
\begin{figure}[h!]
\begin{center}
\subfigure[$\gamma_{1}$]{
\resizebox{1in}{.7in}{\includegraphics{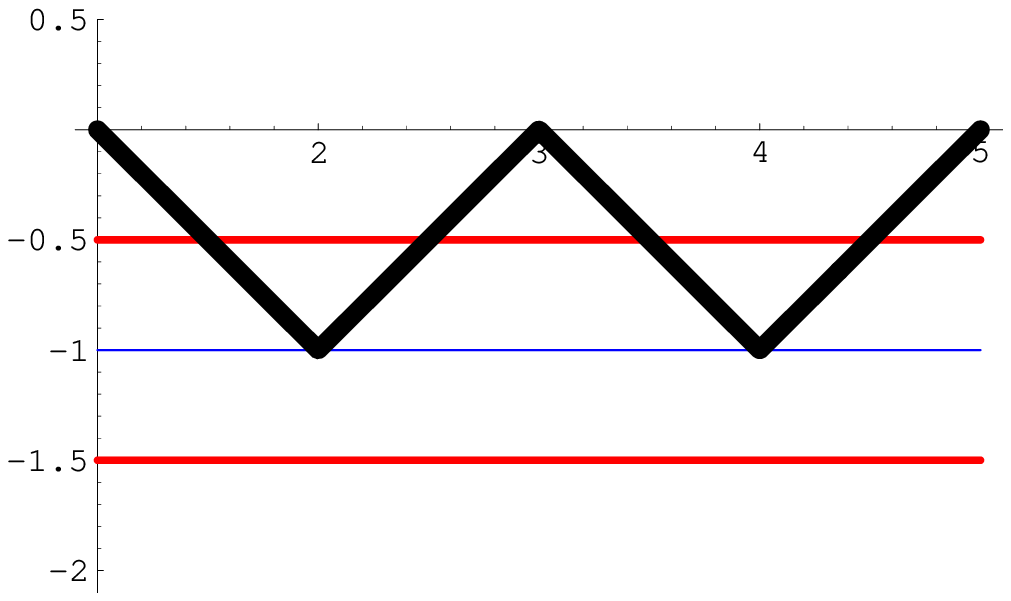}}}
\subfigure[$\gamma_{2}$]{
\resizebox{1in}{.7in}{\includegraphics{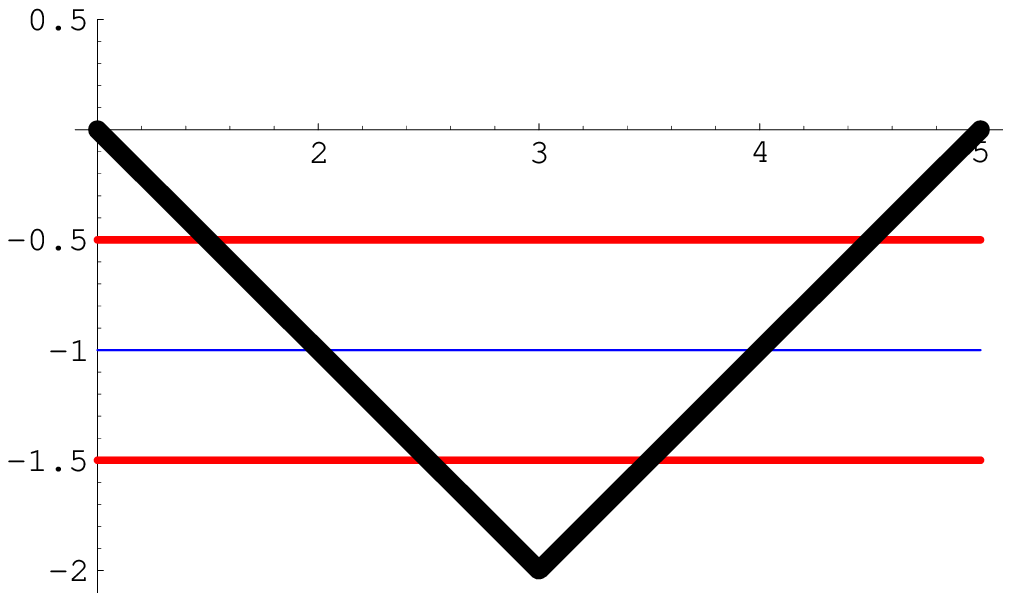}}}
\subfigure[$\gamma_{3}$]{
\resizebox{1in}{.7in}{\includegraphics{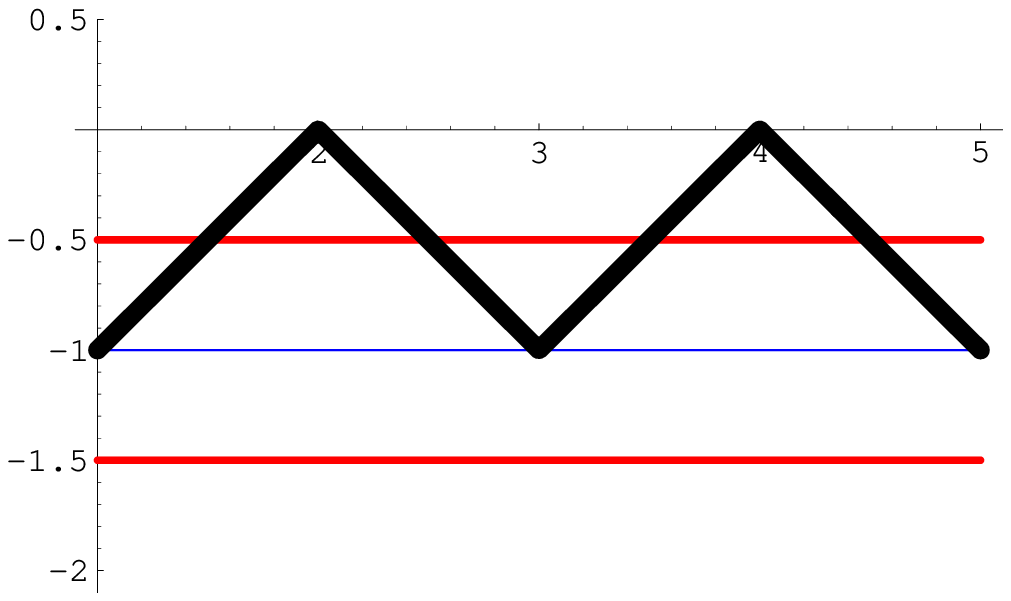}}}
\subfigure[$\gamma_{4}$]{
\resizebox{1in}{.7in}{\includegraphics{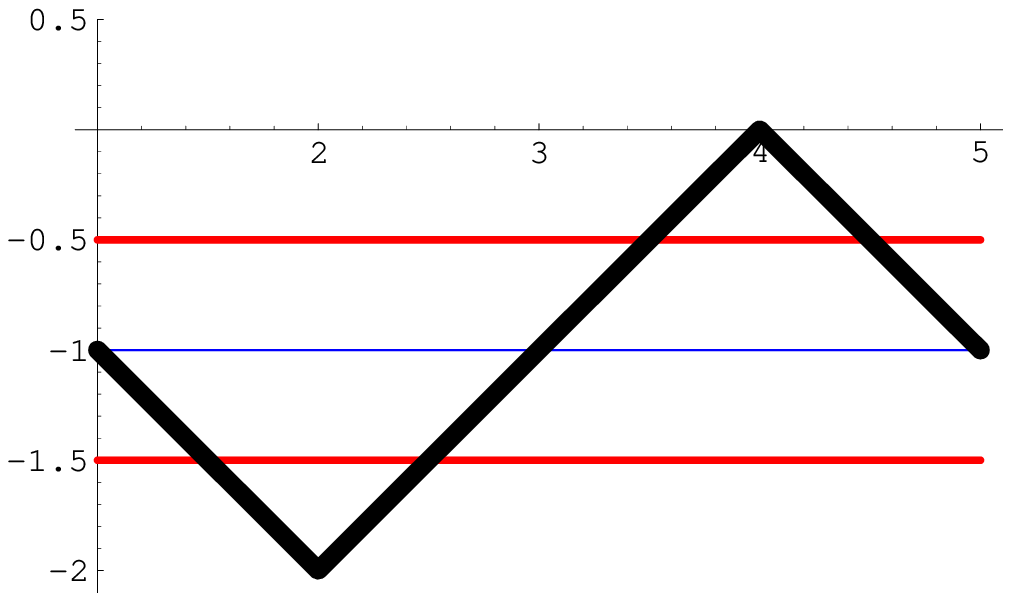}}}
\subfigure[$\gamma_{5}$]{
\resizebox{1in}{.7in}{\includegraphics{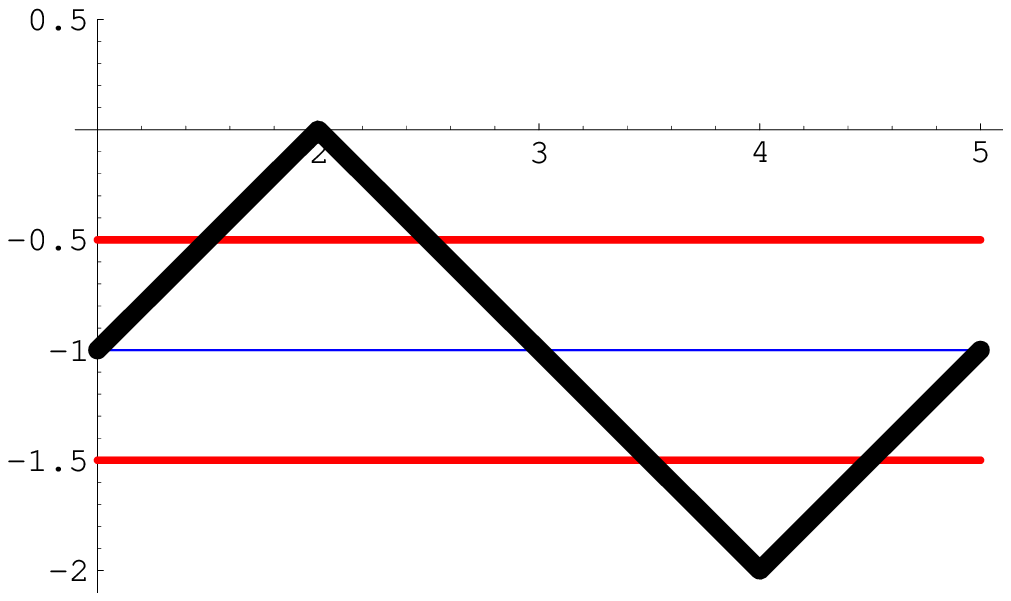}}}
\subfigure[$\gamma_{6}$]{
\resizebox{1in}{.7in}{\includegraphics{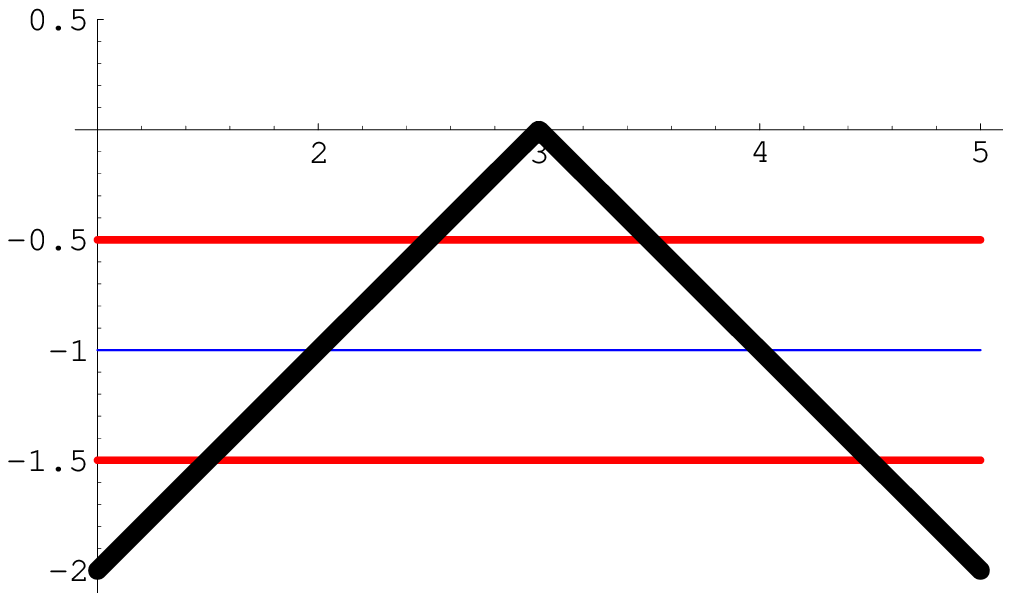}}}
\end{center}
\end{figure}

For each of these strings, one can compute the limit.  For instance, in the case of $\gamma_{1}=(0,-1,0,-1,0)$, according to \eqref{e:1},
\[
\lim_{p\to\infty}\frac{1}{p^{2}}\E[a_{j_{1}+p,j_{2}+p}a_{j_{2}+p,j_{3}+p}a_{j_{3}+p,j_{4}+p}a_{j_{4}+p,j_{1}+p}]=\lim_{p\to\infty}\frac{1}{p^{2}}\E[b_{p-1}^{4}]=1.
\]
Similarly we get $1$ for all the other terms corresponding to these $6$ strings, therefore one gets 
\[
\lim_{n\to\infty}\E[\tr_{n}(X_{n}^{4})]=2,
\] 
which is the fourth moment of the semicircular law $\frac{1}{2\pi}\mathbbm{1}_{[-2,2]}(x)\sqrt{4-x^{2}}dx$.  

For fluctuations, the general statement is give in Theorem~\ref{fluct}.  To be in tune with the convergence discussed thus far, we want to show how one can deal with fluctuations on the following calculation: 
\[
\lim_{n\to\infty}\E\left[ (\Tr_{n}(X_{n}^{4})-\E[\Tr_{n}(X_{n}^{4}))^{2}\right].
\]
Using \ref{eq1},  we can write (neglecting a finite number of terms)
\begin{align*}
 \E&\left[(\Tr_{n}(X_{n}^{4})-\E[\Tr_{n}(X_{n}^{4}))^{2}\right] \sim \frac{1}{n^{4}}\sum_{p,q=1}^{n} S_{p,q},\quad \text{where} \\
 S_{p,q} & = \sum_{\substack{ i_{1},i_{2},i_{3},i_{4} \:\text{admissible},\: \max(i_{1},i_{2},i_{3},i_{4})=p  \\ i'_{1},i'_{2},i'_{3},i'_{4} \:\text{admissible},\: \max(i'_{1},i'_{2},i'_{3},i'_{4})=q } }\E\big[(a_{i_{1},i_{2}}a_{i_{2},i_{3}}a_{i_{3},i_{4}}a_{i_{4},i_{1}}-\E[a_{i_{1},i_{2}}a_{i_{2},i_{3}}a_{i_{3},i_{4}}a_{i_{4},i_{1}}]) \\ & \phantom{aaaaaaaaaaaaaaaaaaaaaaaaaaaaaaaaa} \times (a_{i'_{1},i'_{2}}a_{i'_{2},i'_{3}}a_{i'_{3},i'_{4}}a_{i'_{4},i'_{1}}-\E[a_{i'_{1},i'_{2}}a_{i'_{2},i'_{3}}a_{i'_{3},i'_{4}}a_{i'_{4},i'_{1}}])\big].
\end{align*}
Combining the independence of the entries with the expectation ,we obtain that $S_{p,q}=0$ if $|p-q|\ge2$.  Notice also that $S_{p,q}$  depends only on $p$ and $q$ for $|p-q|\le 1$ and not on $n$.   Therefore, using the following elementary fact 
\[
\lim_{p\to\infty}\frac{x_{p}}{p^{3}}=M \Longrightarrow \lim_{n\to\infty}\frac{1}{n^{4}}\sum_{p=p_{0}}^{n}x_{p}=\frac{M}{4},
\]
it suffices to deal with
\[
\lim_{p\to\infty} \frac{S_{p,q}}{p^{3}}
\]
for $q=p$ or $q=p\pm1$.  We can assume that $p\ge q$, otherwise we can simply reverse the order of $p$ and $q$ for the following argument. Using the same route as before for the computation of the limit, it suffices to find 
\begin{align}\label{f1}
\notag \lim_{p\to\infty} \frac{1}{p^{3}}\big(&\E[a_{j_{1}+p,j_{2}+p}a_{j_{2}+p,j_{3}+p}a_{j_{3}+p,j_{4}+p}a_{j_{4}+p,j_{1}+p}a_{j'_{1}+p,j'_{2}+p}a_{j'_{2}+p,j'_{3}+p}a_{j'_{3}+p,j'_{4}+p}a_{j'_{4}+p,j'_{1}+p}] \\ &\quad-\E[a_{j_{1}+p,j_{2}+p}a_{j_{2}+p,j_{3}+p}a_{j_{3}+p,j_{4}+p}a_{j_{4}+p,j_{1}+p}]\E[a_{j'_{1}+p,j'_{2}+p}a_{j'_{2}+p,j'_{3}+p}a_{j'_{3}+p,j'_{4}+p}a_{j'_{4}+p,j'_{1}+p}]\big), 
\end{align}
where the strings here $j_{1},j_{2},j_{3},j_{4}$ and $j'_{1},j'_{2},j'_{3},j'_{4}$ are admissible with $\max (j_{u},j'_{v},u=1\dots 4,v=1\dots 4)=0$. 

Any appearance of equal consecutive indices in the sequence $j_{1},j_{2},j_{3},j_{4},j_{1}$ or $j'_{1},j'_{2},j'_{3},j'_{4},j'_{1}$, forces another appearance of consecutive indices and such an occurrence means that at least two of the entries in the above limit contains diagonal term.  Since the diagonal terms are bounded (in moments), this forces the above limit to be $0$.  

This implies that the entries in the limit \eqref{f1} are from the subdiagonals only and $(j_{1},j_{2},j_{3},j_{4},j_{1})$ is one of the paths $\gamma_{u}$ in \eqref{gs} and $(j'_{1},j'_{2},j'_{3},j'_{4},j'_{1})$ is one of the paths $\gamma_{u}$ or $\gamma_{u}-1$.  Let's assume for simplicity that we deal with the model in which $b_{n}$ has a $\chi_{n\beta}/\sqrt{\beta}$ distribution.  Using the fact that $\chi_{r}-\sqrt{r}\sim N(0,1/2)$ in distribution and moments sense, one gets that for fixed $k,l\ge1$, $\cov(b_{p}^{k},b_{p}^{l})\sim klp^{(k+l-2)/2}\beta^{-(k+l)/2}/4$. Using this, one can compute the limit above.   For example, if $(j_{1},j_{2},j_{3},j_{4},j_{1})=\gamma_{2}$ and $(j'_{1},j'_{2},j'_{3},j'_{4},j'_{1})=\gamma_{3}-1$, then the limit in \eqref{f1} is 
\begin{align*}
\lim_{p\to\infty}\frac{1}{p^{3}}\left(\E[b_{p}^{2}]\E[b_{p-1}^{6}]-\E[b_{p}^{2}]\E[b_{p-1}^{2}]\E[b_{p-1}^{4}]\right)=\lim_{p\to\infty}\frac{1}{p^{3}}\E[b_{p}^{2}]\cov(b_{p-1}^{4},b_{p-1}^{2})=2\beta^{-3}.
\end{align*}

Using this argument combined with the judicious counting of the paths, one can prove that the fluctuations converge to a Gaussian family.

Let us return  now to the convergence of the empirical distribution of the model \ref{def}.  Since there is nothing sacrosanct about \eqref{e:1}, we can replace it by
\[
\lim_{n\to\infty}\E[(b_{n}/\sqrt{n})^{k}]=m_{k},
\]
where $m_{k}$ is a given number.  Loosely speaking this says that  $b_{n}/\sqrt{n}$ converges in distribution to a random variable $Y$ with moments given by $m_{k}$.  In this case one gets for $\gamma_{1}=(0,-1,0,-1,0)$ that 
\[
\lim_{p\to\infty}\frac{1}{p^{2}}\E[a_{j_{1}+p,j_{2}+p}a_{j_{2}+p,j_{3}+p}a_{j_{3}+p,j_{4}+p}a_{j_{4}+p,j_{1}+p}]=\lim_{p\to\infty}\frac{1}{p^{2}}\E[b_{p-1}^{4}]=m_{4}
\]
and in general, collecting all terms, one gets
\[
\lim_{p\to\infty}\frac{1}{p^{2}}\E[S_{p}]=2m_{4}+4m_{2}^{2}.
\]
The contribution of the paths is as follows: 
\begin{align*}
\gamma_{1}&\to m_{4}, \: & \gamma_{2}&\to m_{2}^{2},\: &\gamma_{3}&\to m_{4},\\ \gamma_{4}&\to m_{2}^{2},\: &\gamma_{5}&\to m_{2}^{2},\:&\gamma_{6}&\to m_{2}^{2}.
\end{align*}
For example $\gamma_{1}$ crosses the line $-1/2$ exactly $4$ times and that corresponds to the index $4$ in $m_{4}$, while the path $\gamma_{2}$ crosses the lines $-1/2$ and $-3/2$ twice, each of these giving an $m_{2}$ term with the total contribution being the product of these, namely $m_{2}^{2}$.

Here we note that the scaling $\sqrt{n}$ in \eqref{e:1} is not essential for the argument.  A more general treatment is one in which $\sqrt{n}$ is replaced by $n^{\alpha}$ with $\alpha>0$, and on this line of ideas the first result we prove is Theorem~\ref{T1} in Section~\ref{1} which concerns the convergence of the traces of powers, both in expectation and almost surely.  There is also a combinatorial relationship between the moments of the limiting distribution and the moments of the limit $b_{n}/n^{\alpha}$ via counting the number of level crossing for paths.  We specialize the limiting distribution in the case $b_{n}/n^{\alpha}$ coverges to $1$.  As opposed to the Wigner ensembles we get here different distributions depending on the scaling used and in some cases even an explicit formula. In Proposition~\ref{ex} we give examples of limiting distributions for the case $b_{n}/n^{\alpha}$ converges to a Bernoulli random variable.  Also worth mentioning here is the fact that the limiting distribution can be described as the distribution in a certain sense of a random Jacobi operator.  At the end  of Section~\ref{1} we also discuss the first order deviation of the expectation of the moments of the distribution of eigenvalues.   

The convergence of the fluctuations is fully discussed in Section~\ref{2}.  Under the appropriate conditions and after properly scaled, the family $\{\tr_{n}(X_{n}^{k})-\E[\tr_{n}(X_{n}^{k})]\}_{k\ge1}$ is shown to converge to a Gaussian family where the covariance can be computed.   

In Section~\ref{3} we extend Theorems~\ref{T1} and \ref{fluct} to the cases of multiple tridiagonal random matrices.  This resembles very much the framework of free probability distribution (see \cite{VDN} for basics and more) and also the second  order freeness discovered by Speicher and Mingo in \cite{SM}.  The interesting part would be to define some kind of cumulant similar to the classical cumulants or to the free cumulant (cf. \cite{SPC}) and then define some sort of ``independence'' via properties of cumulants.  

Section~\ref{4} gives various situations in which the same arguments can be employed to extend  Theorems~\ref{T1} and \ref{fluct}. As a particular case is the band diagonal and an eventual extension to the case in which the entries of the matrix are not independent.  

The combinatorics in this paper is one for Dyke paths.  In the case of Wigner ensembles there is another combinatorial approach by counting planar graphs as it is done in \cite{AZ}.   It would be interesting to see the connection between these two combinatorial methods, though we do not have a clear way of bringing them together.  

The study of tridiagonal models in which dependence of the entries is allowed is very important.  This is motivated for once by studying Wigner ensembles via tridiagonalization.  In this case the independence of the entries of the tridiagonal model is in general lost.  We hope that further study of the tridiagonal models may turn useful in the study of other random matrix models as for example band models where the band width grows with the dimension.   

The general belief is that the tridiagonal random matrices are easier to understand.  This paper is a materialization of this belief in one instance, the case when the entries are independent.   We hope that further study will turn this belief into a scientific fact.

{\bf Acknowledgements}  Special thanks go to Elton Hsu for his interest in this problem and valuable discussions during the preparation of this paper.   Also many thanks to the referees for their valuable comments which led to improvements of this paper.

\section{Convergence of the Distribution of Eigenvalues}\label{1}

Our approach is combinatorial and as such we will deal with the convergence of the distribution of eigenvalues from the moment points of view.  Thus, for a matrix $A$, the $k$th moment of the empirical distribution of eigenvalues is given by the trace of $A^{k}$.  Hence, we reduce the study of the convergence of the moments of the distribution of eigenvalues to the convergence problem of traces of powers of the matrix.  

We start with the following elementary lemma which will be repeatedly used in this paper.

\begin{lemma}\label{l1}
\begin{enumerate}
\item If $x_{n}$ is a sequence of real numbers, then for any $p_{0}\ge1$, and $s>-1$, we have
\begin{equation}\label{e40}
\limsup_{p\to\infty}\frac{|x_{p}|}{p^{s}}\le M \Longrightarrow \limsup_{n\to\infty}\frac{1}{n^{s+1}}\sum_{p=p_{0}}^{n}|x_{p}|\le \frac{M}{s+1}.
\end{equation}

\item If $x_{n}$ is a sequence of real numbers, then for any $p_{0}\ge1$, and $s>-1$, we have
\begin{equation}\label{e4}
\lim_{p\to\infty}\frac{x_{p}}{p^{s}}=M \Longrightarrow \lim_{n\to\infty}\frac{1}{n^{s+1}}\sum_{p=p_{0}}^{n}x_{p}=\frac{M}{s+1}.
\end{equation}

\item If $\{x_{n}\}_{n=1}^{\infty},\{y_{n}\}_{n=1}^{\infty}$ are sequences of real numbers, 
then for any $p_{0}\ge 1$ and $s,t>-1$, we have that 
\begin{equation}\label{e41}
\limsup_{p\to\infty}\frac{|x_{p}|}{p^{s}}\le M,\:\text{and}\:\limsup_{p\to\infty}\frac{|y_{p}|}{p^{t}}\le M' \Longrightarrow \limsup_{n\to\infty}\frac{1}{n^{s+t+2}}\sum_{p_{0}\le p,q \le n}|x_{p}y_{q}|\le \frac{M M'}{(s+1)(t+1)}.
\end{equation}

\item 
If $\{x_{n}\}_{n=1}^{\infty},\{y_{n}\}_{n=1}^{\infty}$ are sequences of real numbers, and $r_{n}$ is a bounded sequence of positive integer numbers, then for all $p_{0},q_{0}$ and  $s,t>0$, we have that 
\begin{equation}\label{e42}
\lim_{p\to\infty}\frac{x_{p}}{p^{s}}=M,\:\text{and}\:\lim_{p\to\infty}\frac{y_{p}}{p^{t}}=M' \Longrightarrow \lim_{n\to\infty}\frac{1}{n^{s+t+2}}\sum_{\substack{p_{0}\le p\le n,\:q_{0}\le q \le n\\ |p-q|\ge r_{n} }}x_{p}y_{q}=\frac{M M'}{(s+1)(t+1)}.
\end{equation}

\end{enumerate}
\end{lemma}

Before we state the first result of this paper we need to introduce some notations.  A path is a string $\lambda=(j_{1},j_{2},\dots,j_{l})$.  A step of $\lambda$ is a pair $(j_{u},j_{u+1})$.  This is called \emph{up} if $j_{u+1}\ge j_{u}+1$,  \emph{down} if $j_{u+1}\le j_{u}-1$ and a \emph{flat} if $j_{u+1}=j_{u}$.  For $k\ge1$, set
\[
\mathcal{P}_{k}=\{ \lambda=(j_{1},j_{2},\dots,j_{k+1})\in\Z^{k+1}:\:j_{1}=j_{k+1},\:|j_{u}-j_{u+1}|\le1,1\le u \le k\},
\]
for the set of paths starting and ending at the same level, and denote by $\mathcal{P}=\cup_{k\ge1}\mathcal{P}_{k}$ the set of all paths starting and ending at the same level.  We call $\lambda$ simply a path and we can realize this as a piecewise path taking the value $j_{u}$ at $u$.  Now for a given integer $p\in\Z$, we define its shift by $p$ units $\lambda+p=(j_{1}+p,j_{2}+p,\dots,j_{k+1}+p)$ and if $\mathcal{R}$ is a set of paths in $\mathcal{P}$, we denote $\mathcal{R}+p=\{\lambda+p:\lambda\in\mathcal{R} \}$.  Given a subset $\Omega$ of $\Z^{2}$ and a  set of numbers $\{a_{i,j}\}_{(i,j)\in \Omega} $ we extend this to $\{a_{i,j} \}_{i,j\in\Z}$ by setting $a_{i,j}=0$ if $(i,j)\in\Z^{2}\setminus \Omega$ and $a_{\lambda}=a_{j_{1},j_{2}}a_{j_{2},j_{3}}\dots a_{j_{k},j_{k+1}}$.
Finally, for a given path $\lambda\in\mathcal{P}$, we set $\max(\lambda)=\max\{j_{u},1\le u\le k\}$, then $\mathbbm{l}_{i}(\lambda)$ to be the number of crosses of the path $\lambda$ with the line $y=i+1/2$ and $\mathbbm{f}_{i}(\lambda)$ the number of flat steps at level $i$, that is the number of pairs $j_{u},j_{u+1}$ appearing in $\lambda$ with $j_{u}=j_{u+1}=i$.  For example, $\lambda=(-2,-2,-3,-2,-2,-1,0,1,1,0,-1,-2,-1,-2)$ has $\mathbbm{l}_{-3}(\lambda)=2$, $\mathbbm{l}_{-2}(\lambda)=4$, $\mathbbm{l}_{-1}(\lambda)=2$, $\mathbbm{l}_{0}(\lambda)=2$, $\mathbbm{f}_{-2}(\lambda)=2$, $\mathbbm{f}_{1}(\lambda)=1$, and the other values of $\mathbbm{l}_{i}(\lambda)$, $\mathbbm{f}_{i}(\lambda)$ are $0$.  Obviously $\mathbbm{l}_{i+p}(\lambda+p)=\mathbbm{l}_{i}(\lambda)$ and similarly 
$\mathbbm{f}_{i+p}(\lambda+p)=\mathbbm{f}_{i}(\lambda)$.

Next, define
{\allowdisplaybreaks
\begin{align*}
\Gamma_{k}&=\{ \gamma=(j_{1}, j_{2}, \dots, j_{k}, j_{k+1})\in\mathcal{P}_{k}:\:  \max(\gamma)=0,\: |j_{u}-j_{u+1}|=1\:\text{for}\:1\le u \le k  \}\\
\Gamma_{k}^{-}&=\{ \gamma=(j_{1}, j_{2}, \dots, j_{k}, j_{k+1})\in\mathcal{P}_{k}\setminus \Gamma_{k}:\:  \max(\gamma)=0\}\: \\
 \Gamma_{k}^{0}&=\{ \gamma=(j_{1},j_{2}, \dots, j_{k}, j_{k+1})\in\mathcal{P}_{k}:\:  j_{1}=j_{k+1}=0,\:\: |j_{u}-j_{u+1}|=1\:\text{for}\:1\le u \le k  \}\\ 
 \Lambda_{k,n}&=\{\lambda=(j_{1},j_{2},\dots,j_{k+1})\in\mathcal{P}_{k}:\: 1\le j_{u}\le n,\: :1\le u \le k \}\\
\Lambda_{n} & =\cup_{k\ge1}\Lambda_{k,n}\\
\Lambda_{k,n}^{p}&=\{ \lambda=(j_{1},j_{2},\dots, j_{k},j_{k+1})\in\Lambda_{k,n}:  \max(\lambda)=p,\: |j_{u}-j_{u+1}|=1,\: \text{for}\:1\le u \le k\} \\
\Lambda_{k,n}^{p,-}&=\{ \lambda=(j_{1},j_{2},\dots, j_{k},j_{k+1})\in\Lambda_{k,n}\setminus\Lambda_{k,n}^{p}: \max(\lambda)=p\}.
\end{align*}
}
Let's point out a couple of simple properties of these sets. All these paths are strings which move at any given step from the previous one by at most one unit and end with the value they started with.   $\Gamma_{k}$ is the set of all paths of length $k$ with only up or down steps, starting and ending at the same level and staying below the $x$-axis, touching it in at least one point.  Similarly,  $\Gamma_{k}^{-}$ is the set of paths of length $k$ with at least one flat step, starting and ending at the same level and staying below the $x$-axis all the time but touching it in at least one point.  $\Lambda_{k,n}$ is the sets of all paths of length $k$ staying above the $x$-axis but below the line $y=n$.  The sets $\Lambda_{k,n}^{p}$ and $\Lambda_{k,n}^{p}$ over $1\le p\le n$ form a partition of the set $\Lambda_{k,n}$.  Notice here an important property which will be exploited below, namely that  $\lambda-\max{\lambda}\in\Gamma_{k}\cup\Gamma_{k}^{-}$ for any $\lambda\in\mathcal{P}_{k}$.  In particular for $p\ge k/2+1$ and any $\lambda\in\Lambda_{k,n}^{p}$ we have that $\lambda-p\in\Gamma_{k}$.  Similarly, for $p\ge k/2+1$, and  any $\lambda\in\Lambda_{k,n}^{p,-}$, one has that $\lambda-p\in\Gamma_{k}^{-}$.  Therefore, if $p\ge k/2+1$,  
\begin{equation}\label{shift}
\Lambda_{k,n}^{p}=\Gamma_{k}+p\quad\text{and}\quad \Lambda_{k,n}^{p,-}=\Gamma_{k}^{-}+p.
\end{equation}
Consequently,  $\Lambda_{k,n}^{p}$ and $\Lambda_{k,n}^{p,-}$ are independent of $n$ for $p\ge k/2+1$.  This simple property turns out to be an important point in proving the next theorem.  
At last, $\Lambda_{n}$ is the collection of all those paths in $\mathcal{P}$ between the lines $y=0$ and $y=n$.  
With these notations, if $A_{n}=\{ a_{i,j}\}_{1\le i,j\le n}$ is the matrix given in \eqref{def}, then, for any path $\lambda\in\Lambda_{n}$, we have 
\begin{equation}\label{ap}
a_{\lambda}=\left(\prod_{j=1}^{n}d_{j}^{\mathbbm{f}_{j}(\lambda)}\right)\left(\prod_{i=1}^{n-1}b_{i}^{\mathbbm{l}_{i}(\lambda)}\right),
\end{equation}
where we use the convention that $0^{0}=1$.

In what follows, for a matrix $X=(x_{ij})_{i,j=1\dots n}$, we denote $\Tr_{n}(X)=\sum_{i=1}^{n}x_{ii}$ and $\tr_{n}=\frac{1}{n}\Tr_{n}$.   

The first result concerns the convergence of the eigenvalue distribution seen at the moment level.  

\begin{theorem}\label{T1}
Let $\alpha>0$.  Assume that all random variables $d_{n}$ and $b_{n}$ are independent and there exists a sequence $\{m_{k}\}_{k\ge0}$, with $m_{0}=1$ so that
\begin{equation}\label{e1}
\lim_{n\to\infty} \E \left[ \left(b_{n}/n^{\alpha}\right)^{k} \right] = m_{k}\quad\text{for any}\quad k\ge0
\end{equation}
and 
\begin{equation}\label{e1'}
\sup_{n\ge1}\E\left[ \left|d_{n}\right|^{k} \right] < \infty\quad\text{for any}\quad k\ge0.
\end{equation}
Denoting $X_{n}=\frac{1}{n^{\alpha}}A_{n}$, we have that 
\begin{equation}
\lim_{n\to\infty} \E\left[ \tr_{n}(X_{n}^{k}) \right]=L_{k}\quad\text{for any}\quad k\ge0,
\end{equation}
and almost surely,
\begin{equation}\label{e2'}
\lim_{n\to\infty}\tr_{n}(X_{n}^{k})=L_{k}\quad\text{for any}\quad k\ge0.
\end{equation}
Moreover, $L_{k}$ is given by 
\begin{equation}\label{e2}
L_{k}=
\begin{cases}
0 & \text{if}\:k\:\text{is odd}\\
\frac{1}{\alpha k+1}\sum_{\gamma\in\Gamma_{k}} \prod_{i< 0} m_{\mathbbm{l}_{i}(\gamma)} & \text{if}\:k\:\text{is even}. 
\end{cases}
\end{equation}

\end{theorem}

\begin{proof}
Notice that, for $k$ and $n$ fixed, the sets $\Lambda^{p}_{k,n}$ and $\Lambda^{p,-}_{k,n}$, $1\le p \le n$ are disjoint and $\cup_{p=1}^{n} \Lambda_{k,n}^{p}\cup\Lambda_{k,n}^{p,-}=\Lambda_{k,n}$.   

As pointed above after the definitions of various $\Gamma$ and $\Lambda$ sets, for  $n\ge p \ge k/2+1$,  $\Lambda^{p,-}_{k,n}=\Gamma_{k}^{-}+p$ and  $\Lambda^{p}_{k,n}=\Gamma_{k}+p$ which implies that 
$\Lambda^{p,-}_{k,n}$ and $\Lambda^{p}_{k,n}$ are independent of $n$.  

Now we denote the elements of the matrix $A_{n}$ by $\{a_{i,j}\}_{1\le i ,j \le n}$ and then write 
\[
\Tr_{n}(A_{n}^{k})=\sum_{1\le i_{1},i_{2},\dots,i_{k}\le n}a_{i_{1},i_{2}}a_{i_{2},i_{3}}\dots a_{i_{k},i_{1}},
\] 
and since $a_{i,j}=0$ for $| i-j |>1$, it follows that 
\begin{equation}\label{e3}
\Tr_{n}(A_{n}^{k})=\frac{1}{n^{\alpha k}}\sum_{\lambda\in \Lambda_{k,n}}a_{\lambda},
\end{equation}
and then 
\[
\tr_{n}(X_{n}^{k})=\frac{1}{n^{\alpha k+1}}\sum_{p=1}^{n}\bigg(\sum_{\lambda\in \Lambda^{p}_{k,n}}a_{\lambda}+\sum_{\lambda\in \Lambda^{p,-}_{k,n}}a_{\lambda}\bigg).
\]

We apply Lemma~\ref{l1} to compute $\lim_{n\to\infty}\E[\tr_{n}(X_{n}^{k})]$.  To this end let's set 
\[
\begin{split}
S^{p}_{n}=\sum_{\lambda\in \Lambda^{p}_{k,n}}a_{\lambda}, & \quad S^{p,-}_{n}=\sum_{\lambda\in \Lambda^{p,-}_{k,n}}a_{\lambda}, \quad 1\le p\le n\\
S^{p}=\sum_{\lambda\in \Gamma_{k}+p}a_{\lambda}, & \quad S^{p,-}=\sum_{\lambda\in \Gamma^{-}_{k}+p}a_{\lambda}, \quad k/2+1\le p.
\end{split}
\]
Since $k$ is fixed and $S^{p,-}_{n}=S^{p,-}$,  $S^{p}_{n}=S^{p}$ for $k/2+1\le p$, combined with the fact that ignoring a finite number of terms does not change the limit of $\tr_{n}(X_{n}^{k})$, we get 
\begin{equation}\label{e5}
\lim_{n\to\infty}\left(
\E [\tr_{n}(X_{n}^{k})] - \frac{1}{n^{\alpha k+1}}\sum_{p=\llbracket k/2\rrbracket+1}^{n}\big( \E[S^{p}]+\E[S^{p,-}] \big)
\right)=0,
\end{equation}
where here $\llbracket s\rrbracket$, stands for the integer part of $s$.    

Here is the key point of the proof.  Invoking \eqref{e4} and \eqref{e5} we reduce the computation of $\lim_{n\to\infty}\E[\tr_{n}(X_{n}^{k})]$ to the computation of 
\[
\lim_{p\to\infty}\frac{1}{p^{\alpha k}}\E[S^{p,-}]\quad\text{and}\quad\lim_{n\to\infty}\frac{1}{p^{\alpha k}}\E[S^{p}].
\]
To do this, first notice that the sums involved in $S^{p}$ and $S^{p,-}$ are finite, therefore everything reduces to computations of the form
\[
\lim_{p\to\infty}\frac{1}{p^{\alpha k}}\E[a_{\gamma+p}]
\]
where $\gamma$ is in $\Gamma_{k}^{-}$ or $\Gamma_{k}$.  For the case $\gamma\in \Gamma_{k}^{-}$,  according to \eqref{ap}, $a_{\gamma+p}=\left(\prod_{j\le0}d_{j+p}^{\mathbbm{f}_{j}(\gamma)}\right)\left(\prod_{i<0}b_{i+p}^{\mathbbm{l}_{i}(\gamma)}\right)$, the products being finite ones.  Thus, using the independence of the entries we get
\[
\E[a_{\gamma+p}]=\left(\prod_{j\le0}\E\left[d_{j+p}^{\mathbbm{f}_{j}(\gamma)}\right]\right)\left(\prod_{i<0}\E\left[b_{i+p}^{\mathbbm{l}_{i}(\gamma)}\right]\right).
\] 
Since $\gamma\in\Gamma_{k}^{-}$, at least one $\mathbbm{f}_{j}(\gamma)$ is $\ge1$. On the other hand $\sum_{j\ge0}\mathbbm{f}_{j}(\gamma)+\sum_{i<0}\mathbbm{l}_{i}(\gamma)=k$, from which one gets that $k-\sum_{i<0}\mathbbm{l}_{i}(\gamma)\ge1$ and
\[
\begin{split}
\left|\frac{1}{p^{k/2}}\E[d_{\gamma+p}]\right| \le \frac{1}{p^{\alpha \left( k-\sum_{i<0}\mathbbm{l}_{i}(\gamma)\right)}}\prod_{j\le0}\left(\E\left[\left|d_{j+p}^{\mathbbm{l}_{j}(\gamma)}\right|\right] \right)\left( \prod_{i<0}\E\left[\left|\left(b_{i+p}/p^{\alpha} \right)^{\mathbbm{l}_{i}(\gamma)}\right| \right]\right),
\end{split}
\]
which together with \eqref{e1} and \eqref{e1'}, yields  that for $\gamma\in\Gamma_{k}^{-}$
\begin{equation}\label{e5.1}
\lim_{p\to\infty}\frac{1}{p^{\alpha k}}\E[a_{\gamma+p}]=0.
\end{equation}
Moreover, since for $k$ odd, $\Gamma_{k}=\emptyset$, this also shows that 
\[
\lim_{n\to\infty}\E [\tr_{n}(X_{n}^{k})]=0,
\]
which is the first part of \eqref{e2}.  If $k$ is even, $p\ge k/2$ and $\gamma\in\Gamma_{k}$, then $a_{\gamma+p}=\prod_{i<0}b_{p+i}^{\mathbbm{l}_{i}(\gamma)}$, and 
\[
\frac{1}{p^{\alpha k}}\E[a_{\gamma+p}]=\prod_{i<0}\E\left[\left( b_{p+i}/p^{\alpha} \right)^{\mathbbm{l}_{i}(\gamma)} \right].
\]
From this, \eqref{e1} and \eqref{e4}, one gets 
\begin{equation}\label{e5.2}
\lim_{p\to\infty}\frac{1}{p^{\alpha k}}\E[a_{\gamma+p}]=\prod_{i<0} m_{\mathbbm{l}_{i}(\gamma)},
\end{equation}
which completes the proof of \eqref{e2}. 

For the almost surely convergence, we use Corollary 1.4.9 from \cite{S1}, which we state here for reader's convenience.  

If $\{X_{n}\}_{n \ge1}$ is a sequence of independent square integrable random variables and  $\{w_{n}\}_{n=1}^{\infty}$ is  a sequence of real numbers which increases to $\infty$ as $n\to\infty$, then, for any $p_{0}\ge 1$, 
\begin{equation}\label{e6}
\sum_{p=1}^{\infty}\frac{\mathrm{var}(X_{p})}{w^{2}_{p}}<\infty  \Longrightarrow \frac{1}{w_{n}}\sum_{p=p_{0}}^{n}(X_{p}-\E[X_{p}])\to 0 \:\text{almost surely}.
\end{equation}
From this, it's  very easy  to deduce the following.  

Let $\{ X_{n}\}_{n\ge1}$ be a sequence of square integrable random variables such that there is an integer constant $q>0$ so that for each $r\in\{0,1,\dots q-1\}$,  $\{X_{r+nq}\}_{n\ge1}$ is a family of independent random variables.  Assume also that $\{w_{n}\}_{n\ge1}$ is a sequence of real numbers which increases to $\infty$ when $n\to\infty$ and has the property that $\lim_{n\to\infty}\frac{w_{n+1}}{w_{n}}=1$. Under these conditions, \eqref{e6} still holds.  
 
In our case, we first point out, that almost surely
\begin{equation}\label{e7}
\lim_{n\to\infty}\left(
\tr_{n}(X_{n}^{k})- \frac{1}{n^{\alpha k+1}}\sum_{p=[k/2]+1}^{n}\big(S^{p}+S^{p,-} \big)
\right)=0.
\end{equation}
Let's notice that for each $r\in\{0,1,\dots,k-1\}$, $\{S^{r+nk}\}_{n\ge1}$ and $\{S^{r+nk,-}\}_{n\ge1}$, are sequences of independent random variables. Now we take $w_{n}=n^{\alpha k+1}$. We show first that 
\[
\sum_{p=k}^{\infty}\frac{\mathrm{var}(S^{p})}{p^{2\alpha k+2}}<\infty,\quad \text{and}\quad \sum_{p=k}^{\infty}\frac{\mathrm{var}(S^{p,-})}{p^{2\alpha k+2}}<\infty,
\] 
which follows once we know that for any $\gamma\in\Gamma_{k}\cup\Gamma^{-}_{k}$, 
\[\tag{*}
\sum_{p=k}^{\infty}\frac{\E\left[a_{\gamma+p}^{2}\right]}{p^{2\alpha k+2}}<\infty.
\] 
To prove this, from \eqref{ap}, $a_{\gamma+p}^{2}=\left( \prod_{j\le 0}a_{j+p}^{2\mathbbm{f}_{j}(\gamma)}\right) \left(\prod_{i<0}b_{i+p}^{2\mathbbm{l}_{i}(\gamma)}\right)$ and using \eqref{e1} and \eqref{e1'} one gets that for a certain constant $C_{k}>0$, 
\[
\frac{1}{p^{2\alpha k}}\E[a_{\gamma+p}^{2}]\le C_{k}\:\text{for}\:p\ge k/2,
\]
which is enough to justify (*) and thus, by \eqref{e6}, that 
\[
 \frac{1}{n^{\alpha k+1}}\sum_{p=[k/2]+1}^{n}\big(S^{p}-\E[S^{p}]+S^{p,-}-\E[S^{p,-}] \big)\xrightarrow[n\to\infty]{}0.
\]
This together with \eqref{e7}, \eqref{e5.1} and \eqref{e5.2} prove \eqref{e2'}. \qedhere
\end{proof}

\begin{remark}
Condition \eqref{e1'} can be relaxed under the assumption that $L_{k}$ are the moments of a measure which satisfy CarlemanÕs 
condition $\sum_{k=0}^{\infty}L_{k}^{-1/(2k)}<\infty$.    In this case, if we replace the condition \eqref{e1'} by the condition that 
\[
\sup_{n\ge1}\E[|d_{n}|]<\infty,
\]
we can conclude that the empirical distribution of eigenvalues of $X_{n}$ converges to the measure whose moments are given by $L_{k}$.  

The proof of this fact is basically given in \cite{B}, page 615 where it is proved that the diagonal part can be removed.  The only essential fact which is needed there is that $\sup_{p\ge1}P(|d_{p}|>n^{\alpha}\epsilon)=o(1)$ which follows from the above condition and Chebyshev's inequality. 
\end{remark}

\begin{corollary}\label{cor1}
Within the notations of the theorem above, assume that $m_{k}=1$ for any even $k$.  Then the numbers $L_{k}$ are the moments of Ullman's distribution  $\nu_{\alpha}(dx)=h_{\alpha}(x)dx$ with ($\alpha>0$)
\[
h_{\alpha}(x)=\mathbbm{1}_{[-2,2]}(x)\frac{1}{\alpha\pi}\int_{|x|/2}^{2}\frac{t^{-1+1/\alpha}}{\sqrt{4-t^{2}}}dt.
\]
These are obtained as distribution of $T^{\alpha}W$, where $W$ has the arsine law distribution ($\mathbbm{1}_{[-2,2]}(x)\frac{1}{\pi\sqrt{4-x^{2}}}dx$) and $T$ is an independent uniform on $[0,1]$.

In some cases, closed formulae are available, as for example,
\begin{equation}\label{exh}
h_{\alpha}(x)=
\begin{cases}
%\mathbbm{1}_{[-2,2]}(x)\frac{1}{\pi \sqrt{4-x^{2}}} & \alpha=0 \\
\mathbbm{1}_{[-2,2]}(x)\frac{ (2+x^{2})\sqrt{4-x^{2}}}{6\pi} & \alpha=1/4 \\
\frac{1}{2\pi}\mathbbm{1}_{[-2,2]}(x) \sqrt{4-x^{2}} & \alpha=1/2 \\
\frac{1}{2\pi}\mathbbm{1}_{[-2,2]}\log\left( (2+\sqrt{4-x^{2}})/|x| \right) & \alpha=1, \\
\end{cases}
\end{equation}
In particular for the model \eqref{e0}, the limiting distribution is the semicircular law. 
\end{corollary}

\begin{figure}[h!]
\begin{center}
\subfigure[$h_{1/4}$]{
\resizebox{1.6in}{1.3in}{\includegraphics{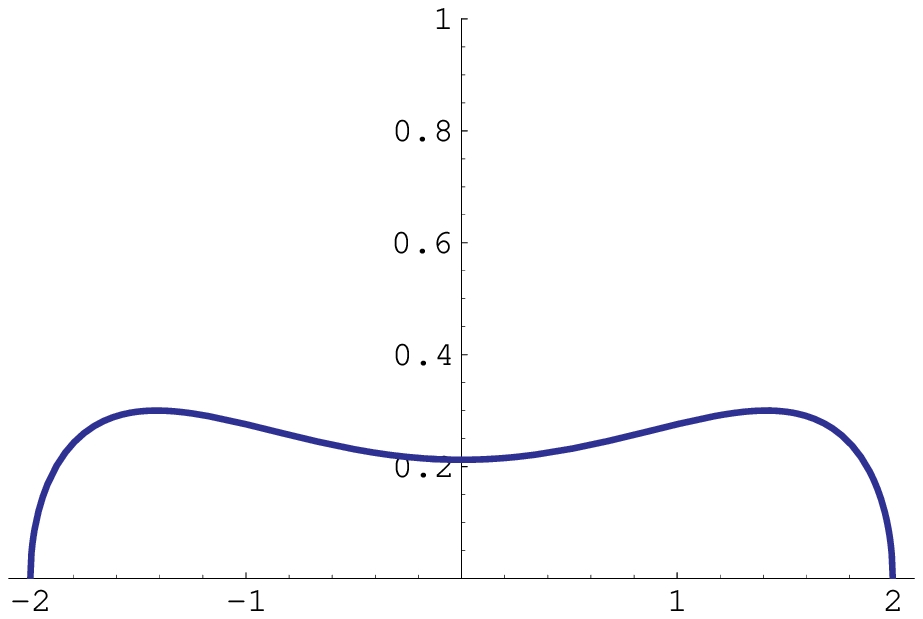}}}
\subfigure[$h_{1/2}$]{
\resizebox{1.6in}{1.3in}{\includegraphics{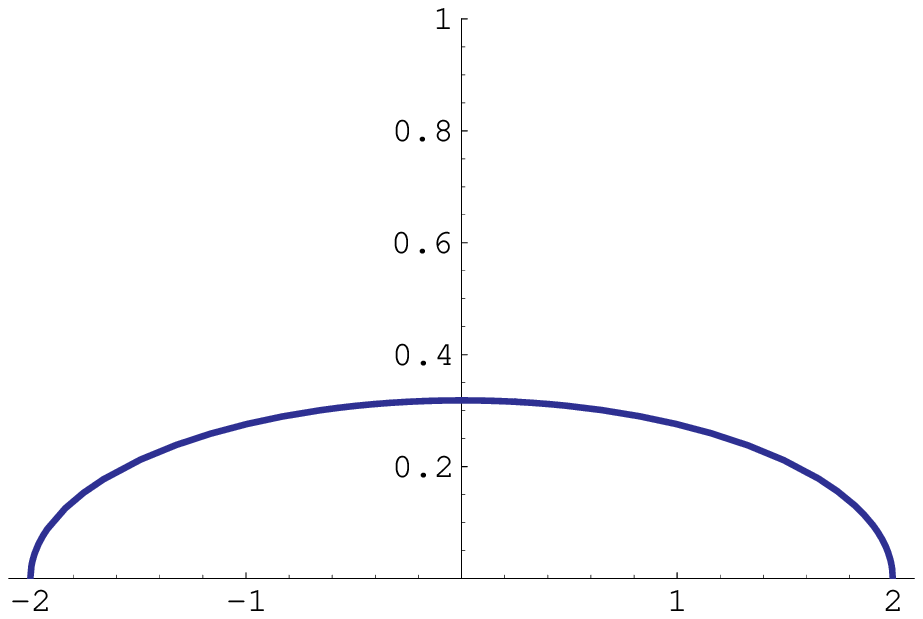}}}
\subfigure[$h_{1}$]{
\resizebox{1.6in}{1.3in}{\includegraphics{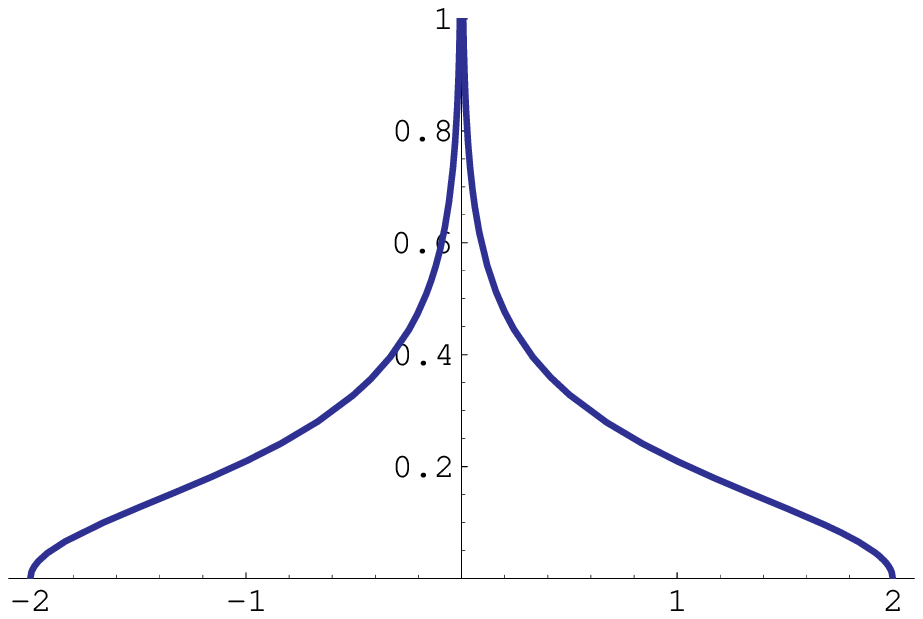}}}
\caption{Density $h_{\alpha}$ for $\alpha=1/4,1/2,1$.}
\end{center}
\end{figure}

\begin{proof}  
 Since $m_{k}=1$ for all $k$ even, we have that the products in \eqref{e2} involving $m$'s equal one.  The number of such terms is given by the number of paths in $\Gamma_{k}$, which turns out to be ${k \choose k/2}$ for $k$ even.  One very quick way to see this is that any path $\gamma$ in $\Gamma_{k}$ is perfectly determined by the prescription of the places where the up steps start, the rest of the positions being filled in with down steps.  Since the path must have the same starting and ending point, it means that there are exactly $k/2$ up steps.  The way of choosing $k/2$ positions out of $k$ points is just 
${k\choose k/2}$.  This means that 
\[\tag{*}
L_{k}=
\begin{cases}
0 & k \: \text{odd}\\
\frac{{k\choose k/2}}{\alpha k+1}&k\: \text{even}.
\end{cases}
\]
For $k$ even and $\alpha=1/2$ these are the moments of the celebrated semicircle law $\frac{1}{2\pi}\mathbbm{1}_{[-2,2]}(x)\sqrt{4-x^{2}}dx$.  

Even though the semicircle plays and important role here, it is the case $\alpha=0$ which is the most important one.  For $\alpha=0$, we have $L_{k}={k\choose k/2}$.  One can check directly that the measure having these properties is the measure $\nu_{0}(dx)=\mathbbm{1}_{[-2,2]}(x)\frac{1}{\pi\sqrt{4-x^{2}}}dx$.  Now if $W$ is a random variable with distribution $\nu_{0}$ and $T$ is an independent and uniform on $[0,1]$, then $T^{\alpha}W$ has the moments given by (*).  From here the rest follows by direct calculations.  
\qedhere 
\end{proof}

\begin{remark}\label{inv}
The system \eqref{e2} is invertible in the sense that for any given the sequence $L_{k}$, $k$ even, one can solve uniquely for the sequence $m_{k}$, $k$ even, since the system \eqref{e2} is a triangular one.  To simplify the notations, set $M_{k}=(\alpha k+1)L_{k}$.  Then we have for the first lines of the system \eqref{e2}: 
\begin{align*}
M_{0}&=m_{0}=1 \\
 M_{2}&=2m_{2} \\
 M_{4}&=2m_{4}+4m_{2}^{2}\\
 M_{6}&=2m_{6}+12m_{4}m_{1}+6m_{2}^{3}\\
 M_{8}&=2m_{8}+16m_{6}m_{2}+12m_{4}^{2}+32m_{4}m_{2}^{2}+8m_{2}^{4}.
 \end{align*}   
 We can solve for $m$'s in terms of $M$'s in this case as:
\begin{align*}
m_{0}&=M_{0}=1 \\
 m_{2}&=\frac{1}{2}M_{2} \\
 m_{4}&=\frac{1}{2}(M_{4}-M_{2}^{2})\\
 m_{6}&=\frac{1}{8}\left(4M_{6}-12M_{4}M_{1}+9M_{2}^{3}\right)\\
 m_{8}&=\frac{1}{4}(2M_{8}-8M_{6}M_{2}-6M_{4}^{2}+28M_{4}M_{2}^{2}-17M_{2}^{4}).
 \end{align*}   
It is of interest a combinatorial interpretation of  this inversion.  Moreover, one such interpretation could perhaps lead to an analytic interpretation, one which would allow characterization of the situation in which the numbers $L_{k}$ are the moments of a real measure.  
 \end{remark}

Next, we would like to compute the limiting distribution in one particular case in which $b_{n}/n^{\alpha}$ converges, not to a constant, but to a Bernoulli random variable. The next proposition also shows that the numbers $L_{k}$ are true moments of a measure under some reasonable conditions.  

Before we state the main result we introduce a class of infinite random matrices known somehow in the theory of random operators as the Anderson model (\cite{T}).  Assume that $\{X_{n}\}_{n\in\Z}$ are given bounded iid random variables.  Then we define 
\begin{equation}\label{op}
\mathcal{A}=\met{\hdotsfor{12}\\ 0&0&X_{-2}&0&X_{-1}&0&0&\hdotsfor{5}\\  \hdotsfor{1}&0&0&X_{-1}&0&X_{0}&0&0&\hdotsfor{4}\\ \hdotsfor{2}&0&0&X_{0}&\fbox{0}&X_{1}&0&0&\hdotsfor{3}\\ \hdotsfor{3}&0&0&X_{1}&0&X_{2}&0&0&\hdotsfor{2}\\  \hdotsfor{4}&0&0&X_{2}&0&X_{3}&0&0&\hdotsfor{1}\\ \hdotsfor{12}}
\end{equation}
where the marked element is the $(0,0)$ element.  We can realize this matrix as a symmetric random Jacobi operator acting on $\ell^{2}(\Z)$. 

Consider the unitary map from $U:\ell^{2}(\Z)\to L^{2}(S^{1})$, where $S^{1}=\{z\in\C:|z|=1\}$ is endowed with the uniform measure and  $U(\{ a_{i}\}_{i\in\Z})(z)=\sum_{i\in\Z}a_{i}z^{i}$. Then the matrix $\mathcal{A}$ becomes the random operator which is given by $\overline{\mathcal{A}}=U\mathcal{A}U^{-1}$ 
\begin{equation}\label{e8}
\overline{\mathcal{A}}(z^{i})=(X_{i}z^{-1}+X_{i+1}z)z^{i}.
\end{equation}
Finally if $e_{i}\in\ell^{2}(\Z)$ is  the vector with $1$ on the $i$th component and $0$ otherwise, then 
\[
A^{k}_{0,0}=\langle \mathcal{A}^{k}e_{0},e_{0} \rangle_{\ell^{2}(\Z)}=\langle \overline{\mathcal{A}}^{k} 1, 1\rangle_{L^{2}(S^{1})},
\]
where here $1$ is the constant function $1$ on $S^{1}$.  

\begin{proposition}\label{ex}
For $\alpha>0$, assume that there is a bounded random variable $Y$ with moments $\E[Y^{2k}]=m_{2k}$ for $k\ge0$.  Then there is a bounded random variable $Z$ whose $k$th moments are $L_{k}$.  

The distribution of $Z$ is of the form $\nu(dx)=\tau f(x)dx+(1-\tau)\delta_{0}$, where $0\le \tau \le 1$ and $f$ is a measurable density function.   

If $\alpha>0$,  $0\le \theta \le 1$ and $m_{2k}=\theta$ for all $k\ge1$, then the distribution whose moments are $L_{k}$, $k\ge0$ is given by
\begin{equation}\label{dist}
\displaystyle{
\nu_{\theta,\alpha}(dx)=
\begin{cases}
\delta_{0}(dx) &  \theta=0 \\
\nu_{\alpha}(dx) & \theta=1\\
 \frac{2\theta}{1+\theta}f_{\theta,\alpha}(x)dx+\frac{1-\theta}{1+\theta}\delta_{0}(dx)& 0<\theta<1,
\end{cases} }
\end{equation}
where 
\[
f_{\theta,\alpha}(x)=\mathbbm{1}_{(-2,2)}(x)\frac{(1+\theta)(1-\theta)^{2}}{2}\sum_{N\ge1}
 \theta^{N-1}g_{N+2}(x),
\]
and for $M\ge 3$
\[
g_{M}(x)=|x|^{1/\alpha-1}\sum_{u=1}^{\llbracket M\arccos(|x|/2)/\pi\rrbracket}\frac{1}{\alpha 2^{1/\alpha}|\cos(u\pi/M)|^{1/\alpha}}.
\]
with the convention that $\sum_{1}^{0}=0$ and $\llbracket s \rrbracket$ stands for the largest integer $\le s$.
\end{proposition}

\begin{figure}[h!]
\begin{center}
\subfigure[$f_{0.4,1/4}$]{
\resizebox{2in}{1.3in}{\includegraphics{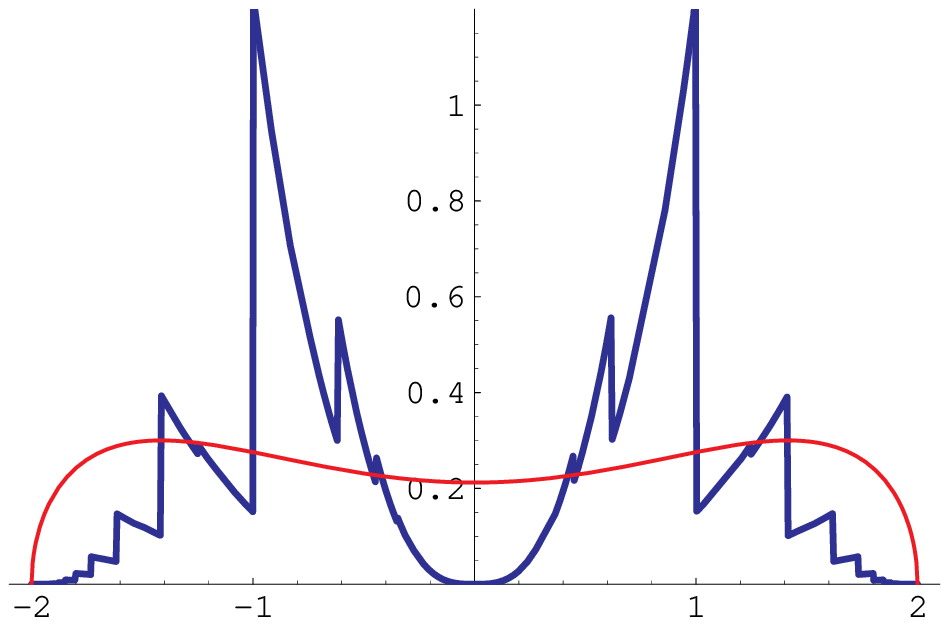}}}
\subfigure[$f_{0.6,1/4}$]{
\resizebox{2in}{1.3in}{\includegraphics{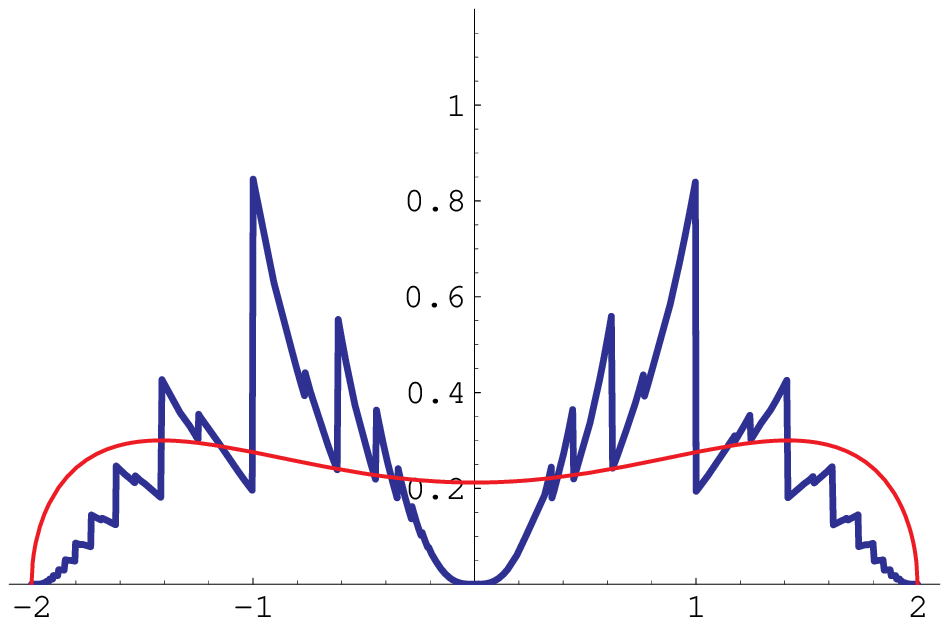}}}
\subfigure[$f_{0.95,1/4}$]{
\resizebox{2in}{1.3in}{\includegraphics{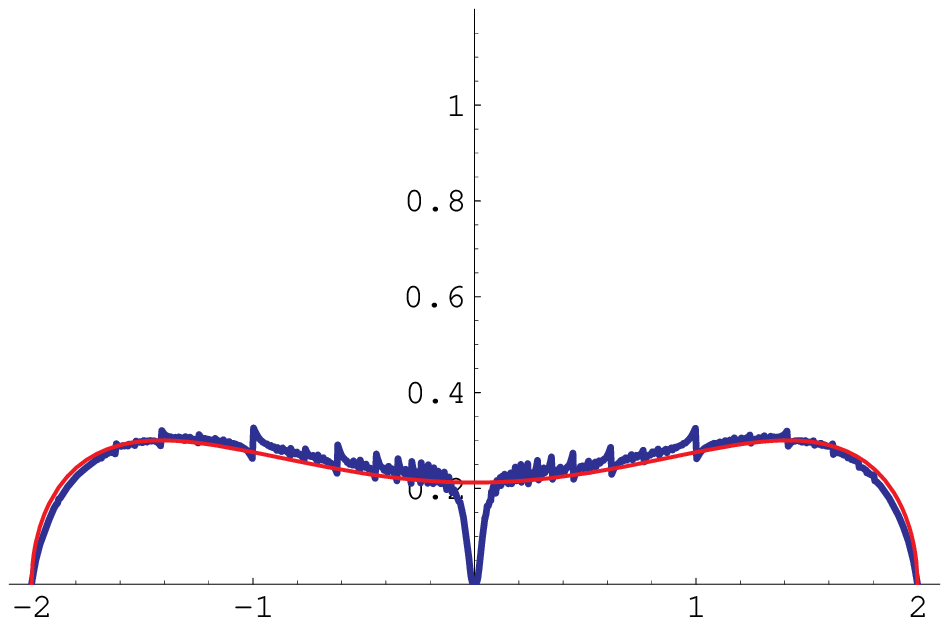}}}
\caption{Density $f_{\theta,1/4}$ for $\theta=0.4,0.6,0.95$ and the smooth density $h_{1/4}$.}

\subfigure[$f_{0.4,1/2}$]{
\resizebox{2in}{1.3in}{\includegraphics{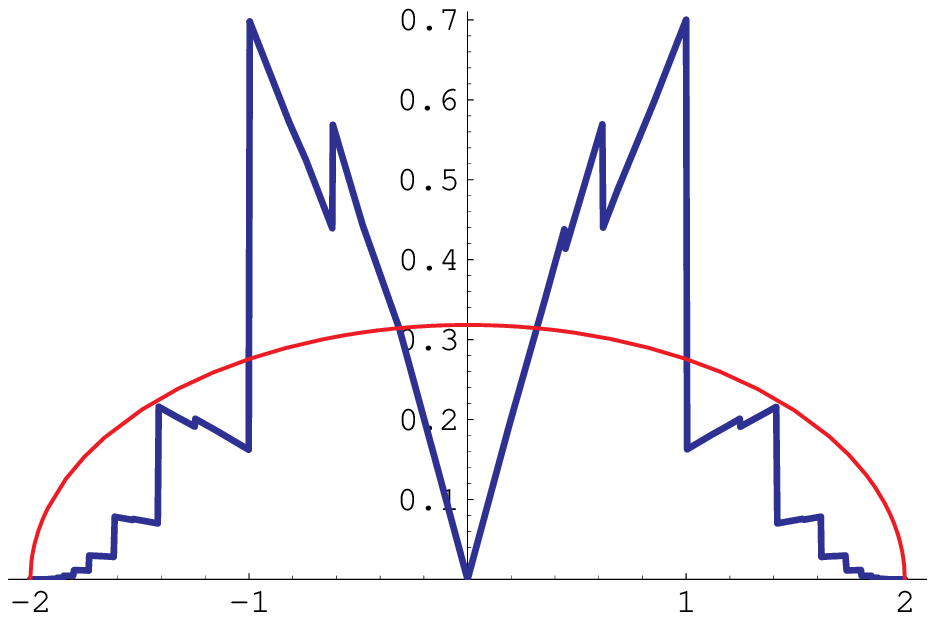}}}
\subfigure[$f_{0.6,1/2}$]{
\resizebox{2in}{1.3in}{\includegraphics{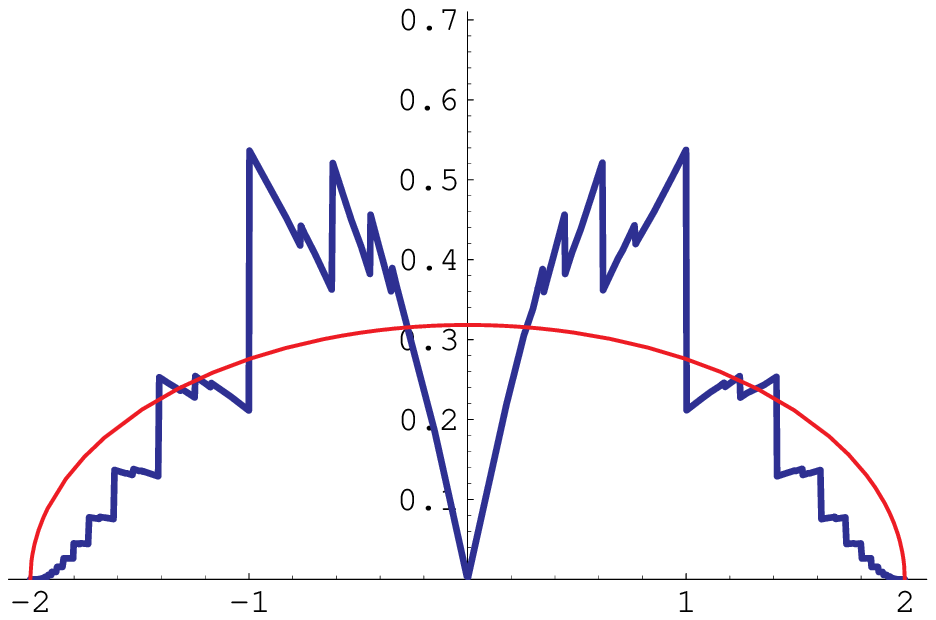}}}
\subfigure[$f_{0.95,1/2}$]{
\resizebox{2in}{1.3in}{\includegraphics{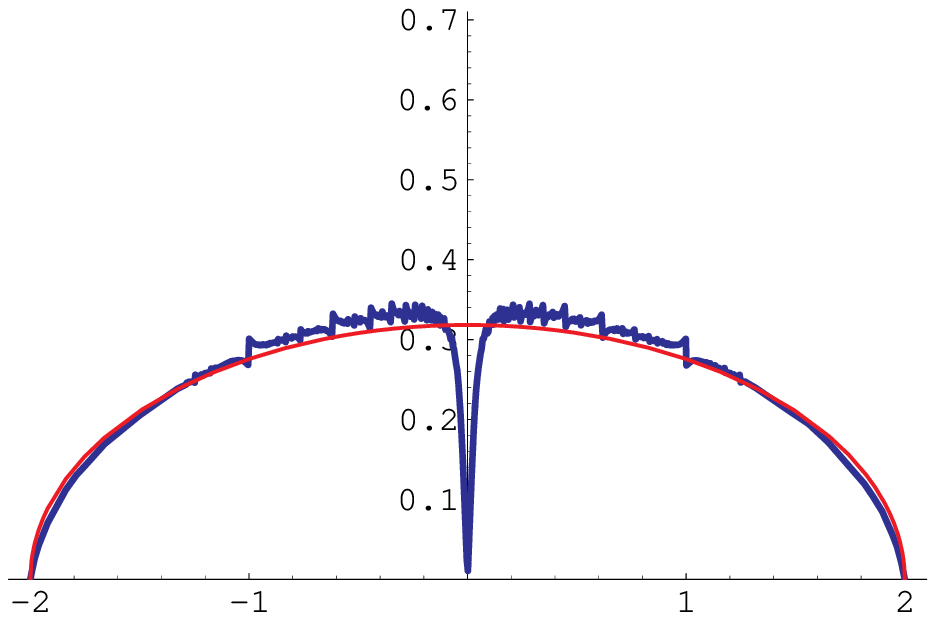}}}
\caption{Density $f_{\theta,1/2}$ for $\theta=0.4,0.6,0.95$ and the smooth semicircular density $h_{1/2}=\frac{1}{2\pi}\sqrt{4-x^{2}}$}

\subfigure[$f_{0.4,1}$]{
\resizebox{2in}{1.3in}{\includegraphics{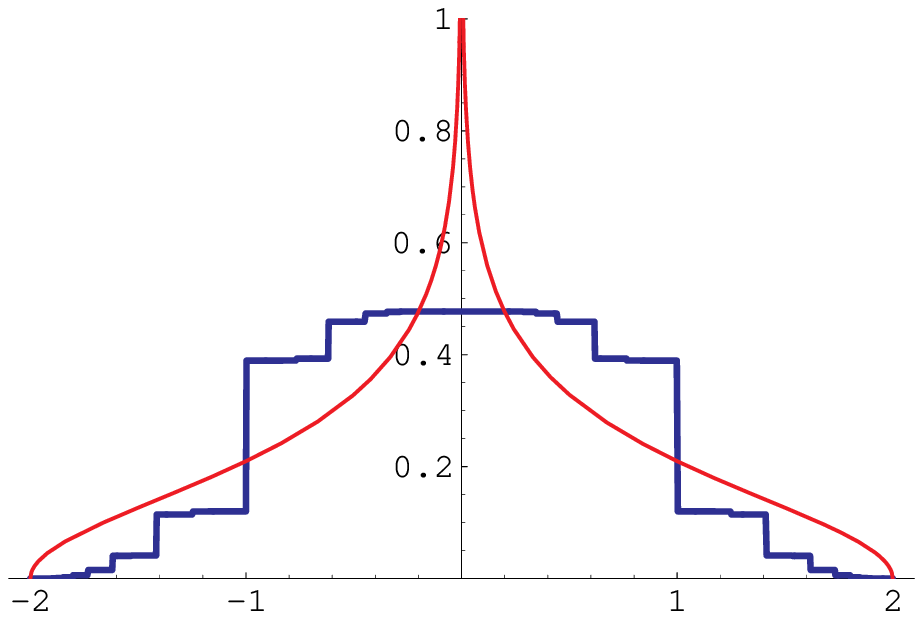}}}
\subfigure[$f_{0.6,1}$]{
\resizebox{2in}{1.3in}{\includegraphics{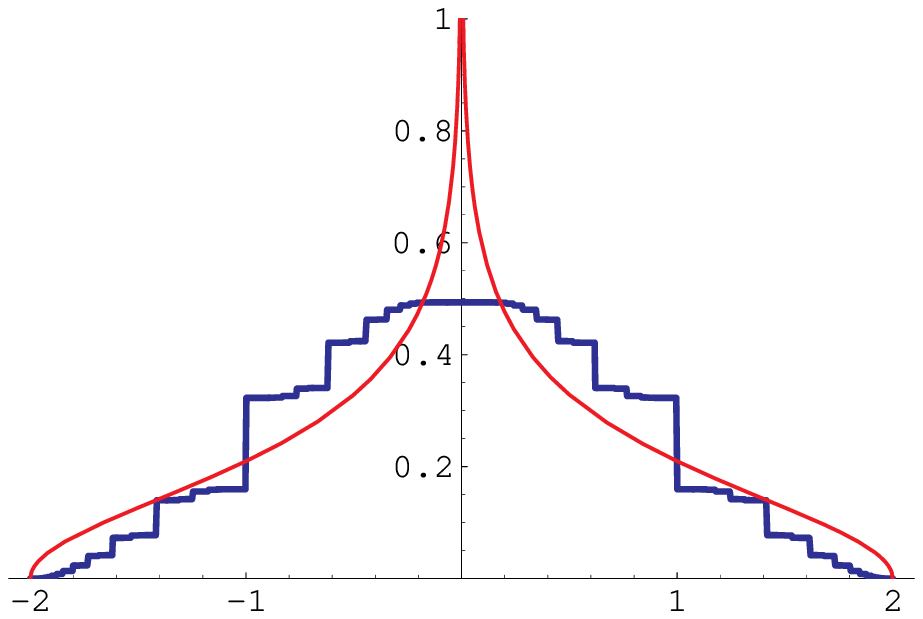}}}
\subfigure[$f_{0.95,1}$]{
\resizebox{2in}{1.3in}{\includegraphics{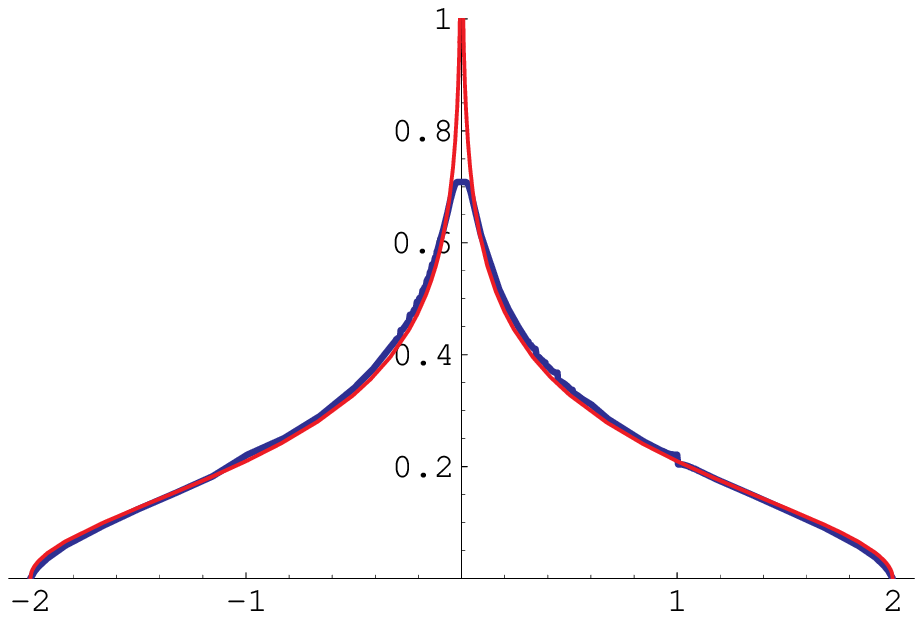}}}
\caption{Density $f_{\theta,1}$, $\theta=0.4,0.6,0.95$ with the smooth density $h_{1}$}

\end{center}
\end{figure}

\begin{proof}
Assume that the distribution of $Y$ is a measure $\mu$ with support in the closed finite interval $I$.  Consider now  the probability space $\Omega=I^{\Z}$ and $P=\mu^{\otimes\Z}$, the product probability on $\Omega$.  We denote by $\omega_{i}$ the $i$th component of $\omega$. Then, define the Hilbert space $H=\{ \mathbf{x}=\{x_{j}\}_{i\in\Z} : x_{j}\in L^{2}(\Omega,P),\sum_{j\in \Z}\| x_{j}\|^{2}_{L^{2}(\Omega,P)}<\infty\}$ with the scalar product given by $\langle\mathbf{x},\mathbf{y} \rangle=\sum_{j\in\Z}\E^{P}[x_{j}y_{j}]$.  On this Hilbert space we consider the operator $\mathbf{A}$ given by 
\[
\mathbf{A}\mathbf{x}=\{ \mathbf{M}_{j}x_{j-1}+\mathbf{M}_{j+1}x_{j+1} \}_{j\in \Z}
\]
for $\mathbf{x}=\{ x_{j}\}_{j\in \Z}$.  Here $\mathbf{M}_{i}: L^{2}(\Omega,P)\to L^{2}(\Omega,P)$ is the multiplication operator given by $(\mathbf{M}_{i}x)(\omega)=\omega_{i} x(\omega)$ for any $x\in L^{2}(\Omega,P)$.  Since $I$ is a closed finite interval, $\mathbf{M}$ is a bounded operator and this in turn yields that the operator $\mathbf{A}$ is also a bounded selfadjoint operator.  

Now we define $\mathbf{e}=\{e_{j} \}_{j\in\Z}$, where $e_{j}=1$ if $j=0$ and $e_{j}=0$ otherwise.  We will prove that 
\begin{equation}\label{A}
\langle \mathbf{A}^{k}\mathbf{e},\mathbf{e} \rangle=
\begin{cases}
0 & k\:\text{odd}\\
\sum_{\gamma\in\Gamma_{k}}\prod_{i<0}m_{l_{i}(\gamma)} &k\:\text{even}
\end{cases}
\end{equation}
To do this we first take the random variables $X_{i}:\Omega\to\R$ given by $X_{i}(\omega)=\omega_{i}$.  The set $\{X_{i} \}_{i\in\Z}$ is a set of iid random variables with distribution $\mu$.  With the random infinite matrix $\mathcal{A}$ given by \eqref{op},  notice that
\begin{equation}\label{eq}
\langle \mathbf{A}^{k}\mathbf{e},\mathbf{e}\rangle=\E[\mathcal{A}^{k}_{0,0}]
\end{equation}
which means that we first compute formally $\mathcal{A}^{k}$ and then take expectation of the $(0,0)$ component.   From this, if we use $a_{i,j}$ for the $(i,j)$ entry of $\mathcal{A}$, then 
\[
\mathcal{A}^{k}_{0,0}=\sum_{i_{1},i_{2},\dots,i_{k}\in\Z}a_{0,i_{1}}a_{i_{1},i_{2}}a_{i_{2},i_{3}}\dots a_{i_{k},0}=
\begin{cases}
0 & k\: \text{odd}\\
\sum_{\gamma\in\Gamma_{k}^{0}}a_{\gamma} & k\:\text{even}.
\end{cases}
\]
This, together with the fact that $a_{i,i+1}=X_{i+1}$, $a_{i-1,i}=X_{i}$, that $\{X_{i}\}_{i\in\Z}$ are iid with distribution $\mu$ and a moment of thinking, gives for $k$ even,
\begin{equation}\label{op1}
\E[\mathcal{A}^{k}_{0,0}]=\sum_{\gamma\in\Gamma_{k}^{0}}\prod_{i\in\Z}\E[X_{i}^{l_{i}(\gamma)}]=\sum_{\gamma\in\Gamma_{k}}\prod_{i<0}m_{l_{i}(\gamma)},
\end{equation}
which proves \eqref{A}.  

On the other hand, since $\mathbf{A}$ is a bounded selfadjoint operator, we can take its spectral measure $\xi(dt)$ and then $\varpi(dt)=\langle \xi(dt)\mathbf{e},\mathbf{e} \rangle$.  We then have that $\int_{\R}t^{k}\varpi(dt)=\langle \mathbf{A}^{k}\mathbf{e},\mathbf{e} \rangle$.  Now if we take a random variable $W$ with distribution $\varpi$ and $T$ an independent uniform random variable on $[0,1]$, one can check, using \eqref{A} and \eqref{e2}, that $Z=T^{\alpha}W$ has the moments $L_{k}$.  It is an easy exercise to verify that the general form of such distributions is $\nu(dx)=\tau f(x)dx+(1-\tau)\delta_{0}$ with $\tau\in[0,1]$.  

For the second part, the case $\theta=0$ is obvious.  Even though the case $\theta=1$ is covered by Corollary \ref{cor1}, we want to employ the arguments used in this proof to reprove it.  The random variable $Y$ in this case is simply the constant $1$.  Therefore, the operator $\mathcal{A}$ becomes a nonrandom operator and, using the representation given by \eqref{e8}, is in fact the multiplication with $z+z^{-1}$ on $L^2({S^{1}})$.  Consequently the spectrum is $[-2,2]$ and 
\[
\int x^{k} \rho(dx)=\langle \overline{\mathcal{A}}^{k} 1, 1\rangle_{L^{2}(S^{1})} = \int_{S^{1}}(z+z^{-1})^{k}dz= \frac{1}{2\pi}\int_{0}^{2\pi}2^{k}(\cos s)^{k}ds=\frac{1}{\pi}\int_{-2}^{2}\frac{x^{k}}{\sqrt{4-x^{2}}}dx,
\]
which results with $\varpi(dx)=\nu_{0}=\frac{1}{\pi}1_{[-2,2]}(x)\frac{1}{\sqrt{4-x^{2}}}$ given in Corrolary~\ref{cor1} .   From this a simple calculation shows that the distribution of $T^{\alpha}W$ where $W$ has distribution $\varpi$ is the $\nu_{\alpha}$ given in Corrolary~\ref{cor1}.  

Next, if $m_{2k}=\theta$ for all $k\ge1$, it is easy to see that $Y$ whose even moments are $m_{2k}$ is a Bernoulli random variable with probability $\theta$ of $1$ and probability $1-\theta$ of $0$.  Thus, the matrix $\mathcal{A}$ has elements $0$ or $1$.  To compute the $(0,0)$ entry of $\mathcal{A}^{k}$, we notice first that (observing when the first $0$ appears in the sequences $\{X_{i}\}_{i\le0}$ and $\{X_{i}\}_{i\ge0}$)
\[
\E[\mathcal{A}^{k}_{0,0}]=\sum_{l,m \ge 0}(\mathcal{A}_{l,m})^{k}_{0,0}\theta^{l+m}(1-\theta)^{2},
\]
where $\mathcal{A}_{l,m}$ is the matrix $\mathcal{A}$, with $X_{0}=X_{-1}=\dots=X_{-m}=1$, $X_{-m-1}=0$ and $X_{1}=X_{2}=\dots=X_{l}=1$, $X_{l+1}=0$.  Thus, the matrix $\mathcal{A}_{l,m}$ is a block matrix of the form 
\[
\mathcal{A}_{l,m}=\met{\mathcal{C}_{l} &\mathbf{0}  & \mathbf{0} \\ \mathbf{0}  &\mathcal{B}_{l,m}&\mathbf{0}  \\ \mathbf{0}  & \mathbf{0} & \mathcal{D}_{m}}
\]
with the square matrix
\[
\mathcal{B}_{l,m}=\met{0&1 &\hdotsfor{5}  \\ 1&0&1& \hdotsfor{4} \\  \hdotsfor{7}\\ \hdotsfor{2} &1&\fbox{0}&1&\hdotsfor{2} \\ \hdotsfor{7} \\ \hdotsfor{4} & 1&0 &1  \\ \hdotsfor{5} & 1 &0 },
\]
 the marked entry being the $(0,0)$ entry of the matrix $\mathcal{A}_{l,m}$ and counting from it, there are $l+1$ rows to the top and $m$ rows to the bottom.  $\mathcal{B}_{l,m}$ is  a $m+l+1$ tridiagonal matrix with only $1$ on the upper and lower sub-diagonals and $0$ otherwise.  The key point is the fact that for $N\ge0$,
\[
\sum_{\substack{l,m\ge0\\ l+m=N}}(\mathcal{A}_{l,m})^{k}_{0,0}\theta^{l+m}(1-\theta)^{2}=\theta^{N}(1-\theta)^{2}\sum_{\substack{l,m\ge0\\ l+m=N}}(\mathcal{B}_{l,m})^{k}_{0,0} =\theta^{N}(1-\theta)^{2}\Tr_{N}\mathcal{B}_{N}^{k},
\]
where $\mathcal{B}_{N}$ is any of the $(N+1)\times (N+1)$ matrices $\mathcal{B}_{l,m}$ with $l+m=N$ and the marked entry removed.  On the other hand,
\[
\Tr_{N}\mathcal{B}_{N}^{k}=\sum_{u=1}^{N}\eta_{u,N}^{k},
\]
where $\eta_{u,N}$ are the eigenvalues counted with their multiplicity.  Now, if $q_{N}(x)=\det(xI_{N}-\mathcal{B}_{N})$ is the characteristic polynomial of $\mathcal{B}_{N}$,  an easy induction argument shows that 
\[
q_{N}(x)=xq_{N-1}(x)-q_{N-2}(x),
\]
with $q_{0}(x)=x$ and $q_{1}(x)=x^{2}-1$.  These, up to scaling, are the Chebyshev polynomials of the second kind.  Precisely, we have $q_{N}(x)=U_{N+1}(\frac{x}{2})$.  As it's well known, the roots of $U_{N}$ are $\cos\left(\frac{u\pi}{N+1}\right)$, $1\le u\le N$ and this shows that the eigenvalues of $\mathcal{B}_{N}$ are 
\[
\eta_{u,N}=2\cos\left(\frac{u\pi}{N+2}\right)
 \: \text{for}\: 1\le u \le N+1.     
\]
From here one gets that 
\[
\E[\mathcal{A}^{k}_{0,0}]=(1-\theta)^{2}\sum_{N\ge0}\theta^{N}\sum_{u=1}^{N+1}\eta_{u,N}^{k}=\int x^{k}\varpi(dx)
\]
where 
\[
\displaystyle
\varpi=\sum_{N\ge 0}\sum_{u=1}^{N+1}(1-\theta)^{2}\theta^{N}\delta_{2\cos\left( \frac{u\pi}{N+2}\right)}.
\]
If this is the distribution of $W$ and $T$ is uniform on $[0,1)$, then $T^{\alpha}W$ has the distribution given by 
\[
\displaystyle{
\nu_{\theta,\alpha}(dx)=\frac{1-\theta}{\theta+1}\delta_{0}(dx)+\sum_{N\ge 1}{\sum_{\substack{u=1\\ u\ne (N+2)/2}}^{N+1}}\frac{|x|^{1/\alpha-1}\theta^{N}(1-\theta)^{2}}{\alpha2^{1/\alpha+1}|\cos(\frac{u\pi}{N+2})|^{1/\alpha}}\mathbbm{1}_{(-2|\cos({u\pi}/{(N+2)})|,2|\cos({u\pi}/{(N+2)})|)}(x)dx }
\]
and from here,  rearrangements bring this to the form given in \eqref{dist}.\qedhere
\end{proof}

Theorem~\ref{T1} gives the zero order convergence in moments of the distribution of eigenvalues.  
Here we are interested in the first order convergence.  The statement can be made more general, but for the sake of simplicity, we give the next theorem in the following form.

\begin{theorem}\label{T2}
Suppose that the conditions \eqref{e1} and \eqref{e1'} from Theorem~\ref{T1} hold.  In addition, assume that there exist $0< \upsilon\le \min(1,2\alpha)$, and numbers $\sigma_{d}$ and $\xi_{k}$ so that 
\begin{equation}\label{e:2}
\lim_{n\to\infty}n^{\upsilon}\left( \E[(b_{n}/n^{\alpha})^{k}]-m_{k} \right)=\xi_{k},\quad k\ge0,
\end{equation}
\begin{equation}\label{e:3}
\E[d_{n}]=0,\:\forall n,\quad\text{and}\quad \lim_{n\to\infty}\var[d_{n}]=\sigma_{d}^{2}
\end{equation}
then, $\lim_{n\to\infty}n^{\upsilon}\left(\E[\tr_{n}(X_{n}^{k})]-L_{k} \right)
$ is given by 
\begin{equation}\label{e:4}
\begin{cases}
\frac{1}{\alpha k -\upsilon + 1}\sum_{\gamma\in\Gamma_{k}}\sum_{j<0}\xi_{\mathbbm{l}_{j}(\gamma)}\prod_{i<0,i\ne j}m_{\mathbbm{l}_{i}(\gamma)}, & k\:\text{even}\:\ge2,\:\upsilon<\min(1,2\alpha)\\
\frac{1}{\alpha k -\upsilon + 1}\big( \sigma_{d}^{2}\sum_{\gamma\in\Gamma_{k}^{2,-}}\prod_{i<0}m_{\mathbbm{l}_{i}(\gamma)}  +\sum_{\gamma\in\Gamma_{k}}\sum_{j<0}\xi_{\mathbbm{l}_{j}(\gamma)}\prod_{i<0,i\ne j}m_{\mathbbm{l}_{i}(\gamma)} \big), & k\:\text{even}\:\ge2,\:\upsilon=2\alpha<1\\
\frac{(\alpha k+1)}{2}L_{k}+\frac{1}{ k}\sum_{\gamma\in\Gamma_{k}}\left(\sum_{i<0}i\mathbbm{l}_{i}(\gamma)\right)\prod_{j<0}m_{\mathbbm{l}_{j}(\gamma)} +\frac{1}{\alpha k }\sum_{\gamma\in\Gamma_{k}}\sum_{j<0}\xi_{\mathbbm{l}_{j}(\gamma)}\prod_{i<0,i\ne j}m_{\mathbbm{l}_{i}(\gamma)},& k\:\text{even}\:\ge2,\:\upsilon=1<2\alpha \\
\frac{( k+2)}{4}L_{k}+\frac{1}{ k}\sum_{\gamma\in\Gamma_{k}}\left(\sum_{j<0}j\mathbbm{l}_{j}(\gamma)\right)\prod_{i<0}m_{\mathbbm{l}_{i}(\gamma)} \\ \qquad+\frac{2}{ k }\sigma_{d}^{2}\sum_{\gamma\in\Gamma_{k}^{2,-}}\prod_{i<0}m_{\mathbbm{l}_{i}(\gamma)}  +\frac{2}{ k}\sum_{\gamma\in\Gamma_{k}}\sum_{j<0}\xi_{\mathbbm{l}_{j}(\gamma)}\prod_{i<0,i\ne j}m_{\mathbbm{l}_{i}(\gamma)}, & k\:\text{even}\:\ge2,\:\upsilon=1=2\alpha\\
0, & k\:\text{odd, or}\:k=0,
\end{cases}
\end{equation}
where $\Gamma_{k}^{2,-}$ is the set of paths in $\Gamma_{k}^{-}$ with exactly two flat steps, both at the same level.

\end{theorem}

\begin{proof}

We use the notations from Theorem~\ref{T1}.  Notice first  that,
\[
\lim_{n\to\infty}n^{\upsilon}\left(
\E [\tr_{n}(X_{n}^{k})] - \frac{1}{n^{\alpha k+1}}\sum_{p=\llbracket k/2\rrbracket+1}^{n}\big( \E[S^{p}]+\E[S^{p,-}] \big)
\right)=0.
\]
Thus, to prove \eqref{e:4} it suffices to find 
\[
\lim_{n\to\infty}\frac{1}{n^{\alpha k-\upsilon+1}}\left(\left(\sum_{p=\llbracket k/2\rrbracket+1}^{n} \E[S^{p}]+\E[S^{p,-}]\right) - n^{\alpha k +1}L_{k}\right),
\]
and, according to Lemma~\ref{l1} and the definition of $L_{k}$,  it reduces to 
\[\tag{*}
\lim_{p\to\infty}\frac{1}{p^{\alpha k-\upsilon}} \E[a_{\gamma+p}]\quad\text{for}\quad\gamma\in\Gamma_{k}^{-} ,
\]
and
\[\tag{**}
\lim_{p\to\infty}\frac{1}{p^{\alpha k-\upsilon}}\left( \E[a_{\gamma+p}] -\frac{(p^{\alpha k+1}-(p-1)^{\alpha k+1})}{\alpha k+1}\prod_{i<0}m_{\mathbbm{l}_{i}(\gamma)}  \right),\quad\text{for}\quad \gamma\in\Gamma_{k}. 
\]
Finally this can be done by using \eqref{ap}. If $\gamma\in\Gamma_{k}^{-}$, then, 
\[
\frac{1}{p^{\alpha k-\upsilon}} \E[a_{\gamma+p}] = \frac{1}{p^{\alpha \sum_{j\le 0}\mathbbm{f}_{j}(\gamma) - \upsilon}}\left(\prod_{j\le 0}\E\left[d_{j+p}^{\mathbbm{f}_{j}(\gamma)} \right]\right)\left(\prod_{i< 0}\E\left[(b_{i+p}/p^{\alpha})^{\mathbbm{l}_{i}(\gamma)} \right]\right).
\]
If $\gamma$ has exactly one flat step, this quantity is zero because $\E[d_{n}]=0$ for all $n$.  On the other hand if $\gamma$ has exactly two flat steps, then for $\upsilon<2\alpha$ this whole term goes to $0$ with $p$.  In the case $\gamma$ has more than $3$ flat steps, one gets that $\sum_{j\le 0}\mathbbm{f}_{j}(\gamma)\ge3$ and so again the term goes to $0$.  The only case we get something nonzero is the case when $\upsilon=2\alpha$ and $\gamma$ has exactly two flat steps at the same level.  In this case, according to \eqref{e:3}, one has that 
\[
\lim_{p\to\infty}\frac{1}{p^{\alpha k-\upsilon}} \E[a_{\gamma+p}]= \sigma_{d}^{2}\prod_{i<0}m_{\mathbbm{l}_{i}(\gamma)}.
\]
Summing over all possible paths in $\Gamma_{k}^{2,-}$,
\[
\lim_{p\to\infty} \frac{1}{p^{\alpha k-\upsilon}}\E[S^{p,-}]=
\begin{cases}
\sigma_{d}^{2}\sum_{\gamma\in\Gamma_{k}^{2,-}}\prod_{i<0}m_{\mathbbm{l}_{i}(\gamma)}, & \upsilon=2\alpha \\
0, & otherwise. 
\end{cases}
\]
If $\gamma\in\Gamma_{k}$, then, 
\begin{align*}
\frac{1}{p^{\alpha k-\upsilon}}\bigg(\E[a_{\gamma+p}] - &\frac{(p^{\alpha k+1}-(p-1)^{\alpha k+1})}{\alpha k+1}\prod_{i<0}m_{\mathbbm{l}_{i}(\gamma)} \bigg) =
\frac{1}{p^{\alpha k  - \upsilon}}\left(\prod_{i< 0}\E\left[b_{i+p}^{\mathbbm{l_{i}}(\gamma)} \right]- p^{\alpha k}\prod_{i<0}m_{\mathbbm{l}_{i}(\gamma)}\right) \\ 
&\quad - \frac{(p^{\alpha k+1}-(p-1)^{\alpha k+1}-(\alpha k+1)p^{\alpha k})}{(\alpha k +1)p^{\alpha k-\upsilon}}\prod_{i<0}m_{\mathbbm{l}_{i}(\gamma)} \\ 
& = \sum_{i<0}\left(\prod_{h< i}\E\left[(b_{h+p}/(p+i)^{\alpha})^{\mathbbm{l}_{h}(\gamma)} \right]\right)
\left(p^{\upsilon}\left(\E\left[(b_{i+p}/p^{\alpha})^{\mathbbm{l_{i}}(\gamma)} \right]- m_{\mathbbm{l}_{i}(\gamma)}\right)\right)\left(\prod_{i<g} m_{\mathbbm{l}_{g}(\gamma)}\right)\\
&\quad+\sum_{i<0}\left(\prod_{h\le i}\E\left[(b_{h+p}/(p+i)^{\alpha})^{\mathbbm{l}_{h}(\gamma)} \right]\right)
\left(p^{\upsilon}\left( (1+i/p)^{\alpha\mathbbm{l}_{i}(\gamma)}-1 \right)\right) \left(\prod_{i\le g} m_{\mathbbm{l}_{g}(\gamma)}\right)\\ &\qquad-\frac{(p^{\alpha k+1}-(p-1)^{\alpha k+1}-(\alpha k+1)p^{\alpha k})}{(\alpha k+1)p^{\alpha k-\upsilon}}\prod_{i<0}m_{\mathbbm{l}_{i}(\gamma)} \\
&\xrightarrow[p\to\infty]{}
\begin{cases}
\sum_{j<0}\xi_{\mathbbm{l}_{j}(\gamma)}\prod_{i<0,i\ne j}m_{\mathbbm{l}_{i}(\gamma)}, &\upsilon<1 \\
\left(\frac{\alpha k}{2}+\alpha\sum_{i<0}i\mathbbm{l}_{i}(\gamma)\right)\prod_{j<0}m_{\mathbbm{l}_{j}(\gamma)}+ \sum_{j<0}\xi_{\mathbbm{l}_{j}(\gamma)}\prod_{i<0,i\ne j}m_{\mathbbm{l}_{i}(\gamma)}, & \upsilon=1.
 \end{cases}
\end{align*}
The rest of \eqref{e:4} follows. \qedhere
\end{proof}

\begin{remark}
We mention here the following equalities for $k$ even $k\ge2$,
\begin{align*} 
|\Gamma_{k}^{2,-}| & =k 2^{k-3},\\ 
\tag{*} \sum_{\gamma\in\Gamma_{k}}\mathbbm{l}_{i}(\gamma)^{2}&=  k 2^{k},\\
\sum_{\gamma\in\Gamma_{k}}\sum_{i<0}i\mathbbm{l}_{i}(\gamma)& =-k2^{k-1}.
\end{align*}
Applying these to the model from \eqref{e0}, we get that $\lim_{n\to\infty}n(\E[\tr_{n}(X_{n}^{k})]-\frac{1}{k/2+1}{k\choose k/2})=(\frac{\beta}{2}-1)(2^{k-1}-{k\choose k/2})$ which are the moments of the measure $\frac{1}{4}(\delta_{-2}+\delta_{2})-\mathbbm{1}_{[-2,2]}(x)\frac{1}{\pi\sqrt{4-x^{2}}}dx$.  This is \cite[Lemma 2.20]{DE} and we just reproved it.  

The proof of the first equation can be done by counting all the paths by directly.  The second equation can be proved using the model $\mathcal{A}(t)$ \eqref{op} with $X_{i}=1+t Y_{i}$ where $Y_{i}$ is a sequence of iid $N(0,1)$ and use \eqref{op1} to get that the identity we are looking at is just the coefficient of $t^{2}$ of $\E[\left(\mathcal{A}(t)\right)^{k}_{00}]$.  Then we can write $\mathcal{A}(t)=\mathcal{B}+t\mathcal{C}$, where $\mathcal{B}$ is the matrix with $1$ on the subdiagonals and $\mathcal{C}$ has iid normal $N(0,1)$ on the subdiagonals. One can compute the coefficient of $t^{2}$ in the $k$th power of $\mathcal{A}(t)$ as a product of the form $\mathcal{B}^{k_{1}}\mathcal{C}\mathcal{B}^{k_{2}}\mathcal{C}\mathcal{B}^{}{k_{3}}$.  The powers $\mathcal{B}^{k}$ can be explicitly computed and then the rest is simple combinatorics.  The proof of the third equality in (*) can be done in the following way.  First realize that the term $\sum_{i<0}i\mathbbm{l}_{i}(\gamma)$ is the negative area between the path and the $x$-axis.  Then one can decompose any path in $\Gamma_{k}$ as two Dyke paths with certain properties.  Finally, one can count the number of paths with a certain area (see \cite[Proposition 6]{F}),  together with manipulations of generating functions to get the equality in there.  

The outlined proofs are long and ad-hoc.  More direct and natural combinatorial proofs are desirable though.          
\end{remark}

\section{Fluctuations}\label{2}

Under the conditions in Theorem \ref{T1} we have almost surely the convergence of the distribution of the eigenvalues of $X_{n}$.  In this section we are interested in the ``fluctuations'' from the limiting distribution.  Theorem \ref{T1}, states that almost surely, 
\[
\lim_{n\to\infty}\left(\tr_{n}X_{n}^{k}-\E[\tr_{n}X_{n}^{k}] \right)=0.
\]
 Next we are interested in how this happens.   More precisely, what is the right factor we should multiply 
 $\tr_{n}X_{n}^{k}-\E[\tr_{n}X_{n}^{k}] $ with to make this converge to something?  Assume that we multiply this by $n^{\eta}$ with $\eta>0$.  What is going to be the right $\eta$?  Let's take a look at the case $k=2$.  Then,
 \[
 n^{\eta}\left(\tr_{n}X_{n}^{2}-\E[\tr_{n}X_{n}^{2}] \right)=\frac{1}{n^{1+2\alpha-\eta}}\sum_{j=1}^{n}\left(d_{j}^{2}-\E[d_{j}^{2}]\right)+\frac{1}{n^{1+2\alpha-\eta}}\sum_{j=1}^{n-1}\left(b_{j}^{2}-\E[b_{j}^{2}]\right).
 \]
Now, for any $0<\eta\le 2\alpha$, the first sum of this goes to $0$ by the Strong Law of Large Numbers.  The other sum can be written as 
\[
\frac{1}{n^{1+2\alpha-\eta}}\sum_{j=1}^{n-1}\left(b_{j}^{2}-\E[b_{j}^{2}]\right)=\frac{1}{n^{1+2\alpha-\eta}}\sum_{j=1}^{n-1}j^{2\alpha}Y_{j}
\]
with $Y_{j}=b_{j}^{2}/j^{2\alpha}-\E[b_{j}^{2}/j^{2\alpha}]$.  

Let's assume that $\{Y_{j}\}_{j=1}^{\infty}$ is a sequence of independent random variables so that in distribution sense $j^{\epsilon}Y_{j}\sim U$ for some $\alpha \ge \epsilon\ge 0$ and $U$ a zero mean random variable with variance $\mathrm{var}(U)>0$.  Then we are looking at the condition that 
\[
Z_{n}=\frac{1}{n^{1+2\alpha-\eta}}\sum_{j=1}^{n-1}j^{2\alpha}Y_{j}\sim \frac{1}{n^{1+2\alpha-\eta}}\sum_{j=1}^{n-1}j^{2\alpha-\epsilon}U_{j} 
\]
is converging in distribution (here $\{U_{j}\}$ are iid with the same distribution as $U$).  Take $\psi(x)$ so that the characteristic function of $U$ is $\E[e^{itU}]=e^{i\psi(t)}$ with $\psi(0)=0$ and $\psi'(0)=0$.  The condition of convergence is translated roughly as convergence when $n\to\infty$ of
\[
\sum_{j=1}^{n-1}\psi(tj^{2\alpha-\epsilon}/n^{1+2\alpha-\eta}).
\]
  Since $U$ is not constant $0$, this implies that $\psi''(0)=\mathrm{var}(U)>0$.  
Now, Taylor expansion $\psi(x)=\psi''(0)x^{2}/2+\mathcal{O}(x^{3})$ yields 
\[
\sum_{j=1}^{n-1}\psi(tj^{2\alpha-\epsilon}/n^{1+2\alpha-\eta}) \sim \mathrm{var}(U)t^{2}\frac{1}{n^{2+4\alpha-2\eta}}\sum_{j=1}^{n-1}j^{4\alpha-2\epsilon}/2\sim
\begin{cases}
\frac{\mathrm{var}(U)t^{2}}{2(1+4\alpha-2\epsilon)} & \eta=1/2+\epsilon \\
0\:\text{or}\:\infty & \text{otherwise}.
\end{cases}
\]
Therefore the choice in this case is obviously $\eta=1/2+\epsilon$.   Moreover, this also shows that the limiting distribution of $Z_{n}$ is normal.  

Another way of guessing $\eta$ is from the general statements of CLT, for variables which are not necessarily identical.

Before we state the next result, we need some definitions. 

We say that the paths $\lambda_{1}, \lambda_{2}\in\mathcal{P}$  \emph{do not share a level} if for any $i\in\Z$, $\mathbbm{l}_{i}(\lambda_{1})\ne0$ implies $\mathbbm{l}_{i}(\lambda_{2})=0$,  $\mathbbm{l}_{i}(\lambda_{2})\ne0$ implies $\mathbbm{l}_{i}(\lambda_{1})=0$ and  similarly for the flat levels, $\mathbbm{f}_{i}(\lambda_{1})\ne0$ implies $\mathbbm{f}_{i}(\lambda_{2})=0$ and $\mathbbm{f}_{i}(\lambda_{2})\ne0$ implies $\mathbbm{f}_{i}(\lambda_{1})=0$.  We say that $\lambda_{1}$ and $\lambda_{2}$ \emph{share a level} if there is an $i$ so that both $\mathbbm{l}_{i}(\lambda_{1})$ are $\mathbbm{l}_{i}(\lambda_{2})$ not zero or both $\mathbbm{f}_{i}(\lambda_{1})$ and $\mathbbm{f}_{i}(\lambda_{2})$ are not zero. 

For $k,l\ge1$, we set $\Gamma(k,l)$ by
\begin{equation}\label{g2}
\begin{cases}
\{ (\gamma_{1},\gamma_{2})\in\mathcal{P}_{k}\times\mathcal{P}_{l}: \max(\max(\gamma_{1}),\max(\gamma_{2}))=0, \gamma_{1},\gamma_{2}\: \text{share a level and have no flat steps} \},& \text{if}\:k,l\:\text{even} \\
\{ (\gamma_{1},\gamma_{2})\in\mathcal{P}_{k}\times\mathcal{P}_{l}: \max(\max(\gamma_{1}),\max(\gamma_{2}))=0, \gamma_{1},\gamma_{2}\: \text{have exactly one flat step each and is shared} \},&\text{if}\: k,l\:\text{odd}.
\end{cases}
\end{equation}

\begin{remark}\label{nrg}
For $k,l$ odd, the number of paths in $\Gamma(k,l)$ is the number of paths of $(\gamma_{1},\gamma_{2})\in\mathcal{P}_{k}\times\mathcal{P}_{l}$, $\gamma_{1},\gamma_{2}$ having exactly one flat step on the $x$-axis.  These pairs can be constructed as follows.  Pick two paths $\gamma_{1}'$ and $\gamma_{2}'$ of length $k-1$ and $l-1$ with only up or down steps.  Then insert any flat step at any level and move the paths so that the level steps are on the $x$-axis.  For $\gamma_{1}'$, there are ${k-1\choose \frac{k-1}{2}}$ choices for the path and $k$ ways of inserting the flat step.  Similarly for $\gamma_{2}'$, so the total number of paths in $\Gamma(k,l)$ is $kl{k-1\choose \frac{k-1}{2}}{l-1\choose \frac{l-1}{2}}$.
\end{remark}

\begin{theorem}\label{fluct}
In addition to  conditions \eqref{e1} and \eqref{e1'} of  Theorem \ref{T1}, assume that,
\begin{equation}\label{dlim}
 \lim_{n\to\infty}\mathrm{var}(d_{n})=\sigma^{2}_{d},
\end{equation}
and there exists $0\le  \epsilon\le \alpha$ so that if $z_{n}^{k}=n^{\epsilon}\left( (b_{n}/n^{\alpha})^{k}-\E[b_{n}/n^{\alpha})^{k}] \right)$, then  
for $k,l\ge 0$ both even,  there exists $C(k,l)$ such that
\begin{equation}\label{ee1'}
\lim_{n\to\infty}\cov\left( z_{n}^{k}, z_{n}^{l}\right) = C(k,l).
\end{equation}
and
\begin{equation}\label{bdd}
\sup_{n\ge1}\E[|z_{n}^{k}|^{m}]<\infty,\forall k,m\ge1.
\end{equation}

Now set
\begin{equation}\label{D}
D(k,l)=
\begin{cases}
\epsilon=0:=
\begin{cases}
\frac{1}{\alpha(k+l)+1}\sum_{(\gamma_{1},\gamma_{2})\in\Gamma(k,l)}\left(\prod_{h<0}m_{\mathbbm{l}_{h}(\gamma_{1})+\mathbbm{l}_{h}(\gamma_{2})}- \prod_{h<0} m_{\mathbbm{l}_{h}(\gamma_{1})}m_{\mathbbm{l}_{h}(\gamma_{2})}\right), &\text{for}\: k,l\:\text{even} \\
0, &\text{otherwise}.
\end{cases} \\
0<\epsilon<\alpha: =
\begin{cases}
\frac{m_{2}^{k+l}}{\alpha(k+l)+1-2\epsilon}\sum_{(\gamma_{1},\gamma_{2})\in\Gamma(k,l)}\sum_{i<0}\frac{C(\mathbbm{l}_{i}(\gamma_{1}),\mathbbm{l}_{i}(\gamma_{2}))}{m_{2}^{\mathbbm{l}_{h}(\gamma_{1})+\mathbbm{l}_{h}(\gamma_{2})}}, &\text{for}\: k,l\:\text{even} \\ 
0, &\text{otherwise}.
\end{cases}  \\
\epsilon=\alpha: =
\begin{cases}
\frac{m_{2}^{k+l}}{\alpha(k+l)}\sum_{(\gamma_{1},\gamma_{2})\in\Gamma(k,l)}\sum_{i<0}\frac{C(\mathbbm{l}_{i}(\gamma_{1}),\mathbbm{l}_{i}(\gamma_{2}))}{m_{2}^{\mathbbm{l}_{h}(\gamma_{1})+\mathbbm{l}_{h}(\gamma_{2})}}, &\text{for}\: k,l\:\text{even} \\
\frac{kl\sigma_{d}^{2}m_{2}^{k+l-2}}{\alpha(k+l)}{k-1 \choose \frac{k-1}{2}}{l-1 \choose \frac{l-1}{2}}, &\: \text{for}\: k,l\:\text{odd},\\ 
0 &\text{otherwise}.
\end{cases}
\end{cases}
\end{equation}
For any polynomial $P(x)=w_{0}+w_{1}x+\dots+w_{N}x^{N}$ 
denote $S_{n}(P)=n^{\epsilon+1/2}\left(\tr_{n}(P(X_{n}))-\E[\tr_{n}(P(X_{n}))]\right)$.  Then 
\[
\lim_{n\to\infty}S_{n}(P)=N(0,\sigma(P)^{2})
\]
where 
\[
\sigma(P)^{2}=\sum_{k,l\ge1}w_{k}w_{l}D(k,l). 
\]
In particular, if $S_{n}(k)=S_{n}(x^{k})$, this implies that the family $\{S_{n}(k)\}_{k\ge1}$ converges in moments to a Gaussian family $\{ S(k)\}_{k\ge1}$ with covariance function $D(k,l)$ and 
\begin{equation}\label{ee2}
\lim_{n\to\infty}S_{n}(k)=
N(0,\sigma_{k}^{2})
\end{equation}
where  $\sigma_{k}^{2}=D(k,k)$.
\end{theorem} 

\begin{remark}\label{r1}
Let's point out that in the case $0<\epsilon\le \alpha$, condition \eqref{ee1'} implies in particular that $m_{k+l}=m_{k}m_{l}$ for any $k,l$ even.  This in turn means that $m_{k}=m_{2}^{k/2}$ for any $k$ even, or that $b_{n}^{2}/n^{2\alpha}\xrightarrow[n\to\infty]{} m_{2}$ in distribution and hence in probability too. 
\end{remark}

\begin{proof}
Write 
\begin{equation}\label{P}
(S_{n}(P))^{j}=\sum_{k_{1},k_{2},\dots,k_{j}=1}^{N}w_{k_{1}}w_{k_{2}}\dots w_{k_{j}}S_{n}(k_{1})S_{n}(k_{2})\dots S_{n}(k_{j}).  
\end{equation}
Since there is a finite number of terms in the above sum, to study the behavior of $\E\left[ \left( S_{n}(P) \right)^{j} \right]$, it suffices to deal with 
\[
\lim_{n\to\infty}\E\left[  S_{n}(k_{1})S_{n}(k_{2})\dots S_{n}(k_{j} )\right]
\]
for a given sequence $k_{1},k_{2},\dots,k_{j}$.  Now, since 
\[
S_{n}(k)=\frac{1}{n^{(\alpha k+(1-2\epsilon)/2)}}\sum_{\lambda\in\Lambda_{k,n}}(a_{\lambda}-\E[a_{\lambda}])
\]
we have 
\begin{equation}\label{f:1}
\E\left[S_{n}(k_{1})S_{n}(k_{2})\dots S_{n}(k_{j} )\right]=\frac{1}{n^{\sum_{i=1}^{j}(\alpha k_{i}+(1-2\epsilon)/2)}}\sum_{\substack{\lambda_{i}\in\Lambda_{k_{i},n} \\ 1\le i \le j}}\E\left[
(a_{\lambda_{1}}-\E[a_{\lambda_{1}}])(a_{\lambda_{2}}-\E[a_{\lambda_{2}}])\dots (a_{\lambda_{j}}-\E[a_{\lambda_{j}}])\right]
\end{equation}
Next we define a notion of connectedness for paths. 
A set of paths $C$ is called \emph{connected} if for any two paths $\lambda$ and $\lambda'$ in $C$ there are paths $\lambda=\lambda_{1},\lambda_{2},\dots,\lambda_{u}=\lambda'$ in $C$ so that  $\lambda_{i},\lambda_{i+1}$ share at least a level.  Otherwise we say that $C$ is not connected or simply disconnected.  The notion of connectedness in this context is an equivalence relation.  Therefore any set $C$ can be written as a disjoint union $C_{1}\cup C_{2}\dots\cup C_{p}$ where each $C_{i}$ is connected.  The sets $C_{1},C_{2},\dots,C_{p}$ are called the connected components of $C$.   If $\lambda$ is a path in $C$, then the connected component containing it is the set of all paths which can be connected with it.  In particular any two paths from different components do not share a level. 

With this concept at hand, we return to \eqref{f:1} and split the sum in sums over all connected components.  Then we organize the connected components in the following way.  For a given partition $\Delta$ of $\{1,2,\dots,j\}$,  we consider $C^{n}_{\Delta}$ the subset of $ (\lambda_{1},\lambda_{2},\dots,\lambda_{j})$ with $\lambda_{i}\in\Lambda_{k_{i},n}$ so that $C^{n}_{\delta}=\{\lambda_{i}:i\in\delta \}$ are the connected components of $C^{n}_{\Delta}$ where $\delta$ runs over all elements of $\Delta$.  In other words the connected components of $(\lambda_{1},\lambda_{2},\dots,\lambda_{j})$ are indexed by the subsets $\delta\in\Delta$.   Now, since any two connected components are disjoint, combined with the independence of the entries of the matrix $A_{n}$, justifies the following rewriting
\[
\E\left[S_{n}(k_{1})S_{n}(k_{2})\dots S_{n}(k_{j} )\right]=\frac{1}{n^{\sum_{i=1}^{j}(\alpha k_{i}+(1-2\epsilon)/2)}}\sum_{\Delta}\sum_{(\lambda_{1},\lambda_{2},\dots,\lambda_{j})\in C_{\Delta}^{n}}\prod_{\delta\in \Delta}\E\left[\prod_{\lambda\in C^{n}_{\delta}}(a_{\lambda}-\E[a_{\lambda}])\right].
\]
Next we fix a partition $\Delta$ of $\{1,2,\dots,j \}$.  The idea is to find the limit of  
\[
U^{\Delta}_{n}=\frac{1}{n^{\sum_{i=1}^{j}(\alpha k_{i}+(1-2\epsilon)/2)}}\sum_{(\lambda_{1},\lambda_{2},\dots,\lambda_{j})\in C_{\Delta}^{n}}\prod_{\delta\in \Delta}\E\left[\prod_{\lambda\in C^{n}_{\delta}}(a_{\lambda}-\E[a_{\lambda}])\right].
\]
To clarify and explain the main idea let's introduce first some notations.   For a given $\mathbf{k}=(k_{1},k_{2},\dots,k_{j})$ and  $\Delta$ a partition of $\{1,2,\dots,j\}$, we set
\[
\Gamma(\mathbf{k}:\Delta)=\{\Gamma=(\gamma_{1},\gamma_{2},\dots,\gamma_{j}): \gamma_{u}\in\mathcal{P}_{k_{u}},\:\text{and for}\:\delta\in\Delta, \max_{u\in\delta}(\max(\gamma_{u}))=0,\: \{\gamma_{u}:u\in\delta\}\:\text{is a connected set} \}.
\]
This is  the set of all paths under the $x$-axis so that by isolating the paths indexed by $\delta$, we obtain a connected set with the maximum of all heights being $0$.  

Notice that for a given $\Delta$, the set $\Gamma(\mathbf{k}:\Delta)$ is actually in a one-to-one correspondence with the set $\times_{\delta\in\Delta}\Gamma(\mathbf{k}_{\delta}:\delta)$, where $\mathbf{k}_{\delta}$ is the vector $\mathbf{k}$ with the components which do not belong to $\delta$ removed.  Obviously there is a finite number of elements in $\Gamma(\mathbf{k}:\Delta)$.  Now if we take a connected component in $C_{\delta}^{n}$, with $\max_{u\in\delta}(\max(\lambda_{u}))=p_{\delta}$, then $\left((\lambda_{u}-p_{\delta})_{ u\in\delta}\right)_{\delta\in\Delta}\in\times_{\delta\in\Delta}\Gamma(\mathbf{k}_{\delta}:\delta)$.  Ignoring eventually a finite number of terms in the expression of $U^{\Delta}_{n}$, the limit of $U^{\Delta}_{n}$ is the same as the limit of 
\[
V^{\Delta}_{n}=\sum_{\Gamma\in \Gamma(\mathbf{k}:\Delta)}\frac{1}{n^{\sum_{\delta\in\Delta}\sum_{u\in\delta}(\alpha k_{u}+(1-2\epsilon)/2)}}\sum_{(p_{\delta})_{\delta\in\Delta}\in\Omega_{n}(\Gamma)}\prod_{\delta\in \Delta}\E\left[\prod_{u\in\delta}(a_{\gamma_{u}+p_{\delta}}-\E[a_{\gamma_{u}+p_{\delta}}])\right]
\] 
where the set $\Omega_{n}(\Gamma)=\{(p_{\delta})_{\delta\in\Delta}:n\ge p_{\delta}\ge p_{0},\:\text{so that}\:\{\gamma_{\delta}+p_{\delta}\}\cap\{ \gamma_{\delta'}+p_{\delta'}\}=\emptyset ,\:\text{for any} \:\delta\ne\delta' \}$, $p_{0}$ being a fixed large number depending only on $\gamma_{\Delta}$ and $\mathbf{k}$.  If the set $\Delta$ has just one element, namely the whole set $\{1,2,\dots,j\}$, then the above sum takes the simpler form
\[
V^{\Delta}_{n}=\sum_{\Gamma\in \Gamma(\mathbf{k}:\Delta)}\frac{1}{n^{\sum_{u\in\{1,2,\dots,j \}}(\alpha k_{u}+(1-2\epsilon)/2)}}\sum_{p=p_{0}}^{n}\E\left[\prod_{u\in\{1,2,\dots,j \}}(a_{\gamma_{u}+p}-\E[a_{\gamma_{u}+p}])\right].  
\]
Using Lemma~\ref{l1}, we can find this limit once one can compute the following
\begin{equation}\label{W}
W(\Gamma)=\lim_{p\to\infty}\frac{1}{p^{\left(\sum_{u\in\{1,2,\dots,j \}}(\alpha k_{u}+(1-2\epsilon)/2)\right)-1}}\E\left[\prod_{u\in\{1,2,\dots,j \}}(a_{\gamma_{u}+p}-\E[a_{\gamma_{u}+p}])\right].  
\end{equation}
Once we know this, we can go back to the case of an arbitrary partition $\Delta$ and use part $3$ and $4$ of Lemma~\ref{l1} to compute the limit of $V^{\Delta}_{n}$. Here are the formulae.  For $\Delta$ with just one component $\{ 1,2,\dots,j\}$, we have
\[
\lim_{n\to\infty}V_{n}^{\Delta}=\sum_{\Gamma\in\Gamma(\mathbf{k}:\Delta)}\frac{2W(\Gamma)}{2\alpha|\mathbf{k}|+j(1-2\epsilon)}
\]
with $|\mathbf{k}|=k_{1}+k_{2}+\dots+k_{j}$.  The general formula which follows from this and a repeated application of Lemma~\ref{l1} is that for an arbitrary partition $\Delta$ we have 
\begin{equation}\label{limU}
\lim_{n\to\infty}U^{\Delta}_{n}=\sum_{\Gamma\in\Gamma(\mathbf{k}:\Delta)}\prod_{\delta\in\Delta}\frac{2W(\gamma_{\delta})}{2\alpha|\mathbf{k}_{\delta}|+|\delta|(1-2\epsilon)}
\end{equation}

Now we want to compute $W(\Gamma)$ when $\Gamma\in\Gamma(\mathbf{k}:\Delta)$ and $\Delta=\{ \{1,2,\dots,j\}\}$.  In the following, for a set $\Omega$ we denote the number of its elements by $|\Omega|$.  

{\bf{Case 1:  $j=1$.}}  In this case due to the fact that  $\E[(a_{\lambda}-\E[a_{\lambda}])]=0$, we get that $U_{n}^{\Delta}=0$ and in particular,$W(\Gamma)=0$ also.

{\bf{Case 2:  $j>2$.}}   We show in this case that $W(\Gamma)=0$.  

To do this we will prove something more general.  Namely we show that for a fixed $\Gamma\in\Gamma(\mathbf{k}:\Delta)$,
\begin{equation}\label{e11}
\E\left[\prod_{u\in\{1,2,\dots,j \}}(a_{\gamma_{u}+p}-\E[a_{\gamma_{u}+p}])\right]=
\mathcal{O}(p^{\sum_{u=1}^{j}(\alpha k_{u}-\epsilon)}) 
\end{equation}

For any path $\gamma$, recall \eqref{ap} which takes the form $a_{\gamma+p}=\left(\prod_{g\le 0 }d_{g+p}^{\mathbbm{f}_{g}(\gamma)}\right)\left(\prod_{h<0}b_{h+p}^{\mathbbm{l}_{h}(\gamma)}\right)$, this product being actually a finite one.  To make the writing in a reasonable form for the expansion of the left hand side in \eqref{e11}, we rewrite 
\[
a_{\gamma+p}=\prod_{i\in\Z }c(i)_{p-|i|}^{\mathbbm{m}_{i}(\gamma)}
\]
where 
\[
c(i)=
\begin{cases}
d & i\le 0 \\
b & i>0
\end{cases},
\quad\text{and}\quad
\mathbbm{m}_{i}= 
\begin{cases}
\mathbbm{f}_{i} & i\le 0 \\
\mathbbm{l}_{-i} & i>0.
\end{cases}
\]
Since the entries are independent,  we have
\[
\begin{split}
a_{\gamma+p}-\E[a_{\gamma+p}] & = \sum_{i\le0}\left(\prod_{g\le i-1 }d_{g+p}^{\mathbbm{f}_{g}(\gamma)}\right)\left(d_{g+i}^{\mathbbm{f}_{i}(\gamma)}-\E\left[d_{g+i}^{\mathbbm{f}_{i}(\gamma)}\right]\right) \left(\prod_{ i+1\ge g }\E\left[d_{g+p}^{\mathbbm{f}_{g}(\gamma)}\right]\right)\left(\prod_{h<0}\E\left[b_{h+p}^{\mathbbm{l}_{h}(\gamma)}\right]\right) \\ & \quad +  
\sum_{i<0}\left(\prod_{g\le 0 }d_{g+p}^{\mathbbm{f}_{g}(\gamma)}\right)\left(\prod_{h<i}b_{h+p}^{\mathbbm{l}_{h}(\gamma)}\right)\left( b_{i+p}^{\mathbbm{l}_{i}(\gamma)}-\E\left[ b_{i+p}^{\mathbbm{l}_{i}(\gamma)}\right] \right)\left( \prod_{i+1\le h}\E\left[b_{h+p}^{\mathbbm{l}_{h}(\gamma)}\right]\right)
\end{split}
\]
which can be rewritten as
\[
a_{\gamma+p}-\E[a_{\gamma+p}]=\sum_{l\in\Z}\prod_{i\in\Z}\left( c(i)_{p-|i|}^{\mathbbm{m}_{i}(\gamma)}\right)^{\tau_{l}(i)}\left( c(i)_{p-|i|}^{\mathbbm{m}_{i}(\gamma)}-\E\left[c(i)_{p-|i|}^{\mathbbm{m}_{i}(\gamma)}\right] \right)^{\nu_{l}(i)}\left(\E\left[c(i)_{p-|i|}^{\mathbbm{m}_{i}(\gamma)}\right] \right)^{\zeta_{l}(i)}
\]
with 
\[
\tau_{l}(i)=
\begin{cases}
1 & i<l \\ 
0 & i\ge l
\end{cases},
\quad
\nu_{l}(i)=
\begin{cases}
1 & i=l \\ 
0 & i\ne l
\end{cases},
\quad
\zeta_{l}(i)=
\begin{cases}
1 & i>l \\ 
0 & i\le l.
\end{cases}
\]
Notice that $\tau_{l}+\nu_{l}+\zeta_{l}=1$, $\sum_{i\in\Z}\nu_{l}(i)=1$ for any $l\in\Z$ and if $\gamma\in\mathcal{P}_{k}$, then $\sum_{i\in\Z}\mathbbm{m}_{i}(\gamma)=k$.  
Using these formulae for $\gamma_{1},\gamma_{2},\dots,\gamma_{j}$, after multiplying out the factors, the left hand side in \eqref{e11} becomes
\[
\sum_{l_{1},l_{2},\dots,l_{j}\in\Z}\prod_{i\in\Z}\E\left[\prod_{u=1}^{j}\left( c(i)_{p-|i|}^{\mathbbm{m}_{i}(\gamma_{u})}\right)^{\tau_{l_{u}}(i)}\left( c(i)_{p-|i|}^{\mathbbm{m}_{i}(\gamma_{u})}-\E\left[c(i)_{p-|i|}^{\mathbbm{m}_{i}(\gamma_{u})}\right] \right)^{\nu_{l_{u}}(i)} \left(\E\left[c(i)_{p-|i|}^{\mathbbm{m}_{i}(\gamma_{u})}\right] \right)^{\zeta_{l_{u}}(i)}\right].
\]
Now, if $i\le 0$, then $c(i)=d$ and then \eqref{e1'} combined with H\" older's inequality yields that each product in the above sum with $i\le 0$ is bounded by a constant.  If $i>0$, then $c(i)=b$, $\mathbbm{m}_{i}=\mathbbm{l}_{-i}$ and in this case, using H\" older's inequality, \eqref{e1}, \eqref{bdd} and \eqref{ee1'}  one can show that 
\[
\E\left[\prod_{u=1}^{j}\left( c(i)_{p-|i|}^{\mathbbm{m}_{i}(\gamma_{u})}\right)^{\tau_{l_{u}}(i)}\left( c(i)_{p-|i|}^{\mathbbm{m}_{i}(\gamma_{u})}-\E\left[c(i)_{p-|i|}^{\mathbbm{m}_{i}(\gamma_{u})}\right] \right)^{\nu_{l_{u}}(i)} \left(\E\left[c(i)_{p-|i|}^{\mathbbm{m}_{i}(\gamma_{u})}\right] \right)^{\zeta_{l_{u}}(i)}\right]
\le C p^{\sum_{u=1}^{j}\left(\alpha\mathbbm{m}_{i}(\gamma_{u})-\epsilon\nu_{l_{u}}(i)\right)},
\]
where the constant $C$ depends only on the paths $\gamma_{1},\gamma_{2},\dots,\gamma_{j}$.   This means that for fixed $l_{1},l_{2},\dots,l_{j}\in\Z$
\[
\prod_{i\in\Z}\E\left[\prod_{u=1}^{j}\left( c(i)_{p-|i|}^{\mathbbm{m}_{i}(\gamma_{u})}\right)^{\tau_{l_{u}}(i)}\left( c(i)_{p-|i|}^{\mathbbm{m}_{i}(\gamma_{u})}-\E\left[c(i)_{p-|i|}^{\mathbbm{m}_{i}(\gamma_{u})}\right] \right)^{\nu_{l_{u}}(i)} \left(\E\left[c(i)_{p-|i|}^{\mathbbm{m}_{i}(\gamma_{u})}\right] \right)^{\zeta_{l_{u}}(i)}\right]\le C p^{\sum_{u=1}^{j}\sum_{i>0}\left(\alpha\mathbbm{m}_{i}(\gamma_{u})-\epsilon\nu_{l_{u}}(i)\right)}.
\]
Next we have 
\[
\sum_{i>0}\left(\alpha\mathbbm{m}_{i}(\gamma_{u})-\epsilon\nu_{l_{u}}(i)\right) \le 
\begin{cases}
\alpha k_{u}-\epsilon & \text{if}\:\gamma_{u}\:\:\text{contains no flat step} \\
\alpha k_{u}-\alpha & \text{if}\:\gamma_{u}\:\:\text{contains exactly one flat step} \\
\alpha k_{u}-2\alpha & \text{if}\:\gamma_{u}\:\:\text{contains two or more flat steps}.
\end{cases}
\]
To see this, one should notice that if $\gamma_{u}$ does not contain a flat step, then $\sum_{i>0}\mathbbm{m}_{i}(\gamma_{u})=\sum_{i<0}\mathbbm{l}_{i}(\gamma_{u})=k_{u}$, while $\sum_{i>0}\nu_{l_{u}}(i)=1$.  In the case $\gamma_{u}$ has just one flat step then, if $l_{u}\le0$, then $\sum_{i>0}\mathbbm{l}_{i}(\gamma_{u})= k_{u}-1$ while $\sum_{i>0}\nu_{l_{u}}(i)=0$ and if $l_{u}>0$, then 
$\sum_{i>0}\mathbbm{l}_{i}(\gamma_{u})= k_{u}-1$ while $\sum_{i>0}\nu_{l_{u}}(i)=1$ which justifies the first part.  In the case $\gamma_{u}$ has more than one flat step, then $\sum_{i>0}\mathbbm{l}_{i}(\gamma_{u})\le k_{u}-2$ and the rest follows.  Hence, since $\epsilon\le\alpha$, 
\begin{equation}\label{le}
\E\left[\prod_{u\in\{1,2,\dots,j \}}(a_{\gamma_{u}+p}-\E[a_{\gamma_{u}+p}])\right] \le 
\begin{cases}
C p^{\sum_{u=1}^{j}(\alpha k_{u}-\epsilon)} & \text{if any}\:\gamma_{u},1\le u \le j,\:\text{has no flat steps} \\
C p^{\left(\sum_{u=1}^{j}(\alpha k_{u}-\epsilon)\right)-(\alpha-\epsilon)} & \text{if any}\:\gamma_{u},1\le u \le j,\:\text{has at most one flat step} \\
&\text{and at least one has exactly one flat step}\\
C p^{\left(\sum_{u=1}^{j}(\alpha k_{u}-\epsilon)\right)-\alpha} & \text{if one of}\: \gamma_{u},1\le u \le j,\:\text{has two or more flat steps}. 
\end{cases}
\end{equation}
which suffices to prove \eqref{e11}.  

{\bf{Case 3:  $j=2$.}}  In this case, $\Gamma=(\gamma_{1},\gamma_{2})$ and we need to compute
\begin{equation}\label{e12}
\lim_{p\to\infty}\frac{1}{p^{\alpha k_{1}+\alpha k_{2}-2\epsilon}}\E[(a_{\gamma_{1}+p}-\E[a_{\gamma_{1}+p}])(a_{\gamma_{2}+p}-\E[a_{\gamma_{2}+p}])].
\end{equation}

Here we distinguish the cases $\epsilon<\alpha$ and $\epsilon=\alpha$.  From equation \eqref{le}, we see that for $\epsilon<\alpha$, the dominant term is the one involving only sums over the paths with no flat steps.  If $\epsilon=\alpha$, then we need to consider also the paths with exactly one flat step.  

First we consider the contribution from the paths with no flat steps.  
To carry this out, invoke \eqref{ap} and since there are no flat steps, 
$a_{\gamma+p}=\prod_{h<0}b_{h+p}^{\mathbbm{l}_{h}(\gamma)}$ and 
\[
a_{\gamma+p}-\E[a_{\gamma+p}] =
\sum_{i<0}\left(\prod_{h<i}b_{h+p}^{\mathbbm{l}_{h}(\gamma)}\right)\left( b_{i+p}^{\mathbbm{l}_{i}(\gamma)}-\E\left[ b_{i+p}^{\mathbbm{l}_{i}(\gamma)}\right] \right)\left(\prod_{i+1\le h}\E\left[b_{h+p}^{\mathbbm{l}_{h}(\gamma)}\right]\right),
\]
from which one gets 
{\allowdisplaybreaks
\begin{align*}
(a_{\gamma_{1}+p}-&\E[a_{\gamma_{1}+p}])(a_{\gamma_{2}+p}-\E[a_{\gamma_{2}+p}])=
\sum_{i_{1}<i_{2}<0}\left(\prod_{h<i_{1}} b_{h+p}^{\mathbbm{l}_{h}(\gamma_{1})+\mathbbm{l}_{h}(\gamma_{2})}\right)  \left(\left( b_{i_{1}+p}^{\mathbbm{l}_{i_{1}}(\gamma_{1})}-\E\left[ b_{i_{1}+p}^{\mathbbm{l}_{i_{1}}(\gamma_{1})}\right] \right)b_{i_{1}+p}^{\mathbbm{l}_{i_{1}}(\gamma_{2})} \right)\\*
&\qquad \times\left(\prod_{i_{1}<h< i_{2}} \E\left[b_{h+p}^{\mathbbm{l}_{h}(\gamma_{1})}\right] b_{h+p}^{\mathbbm{l}_{h}(\gamma_{2})}\right) \left( b_{i_{2}+p}^{\mathbbm{l}_{i_{2}}(\gamma_{2})}-\E\left[ b_{i_{2}+p}^{\mathbbm{l}_{i_{2}}(\gamma_{2})} \right] \right) \left(\prod_{i_{2}<h} \E\left[b_{h+p}^{\mathbbm{l}_{h}(\gamma_{2})}\right]\E\left[ b_{h+p}^{\mathbbm{l}_{h}(\gamma_{1})}\right] \right) \\ 
& + \sum_{i<0}\left(\prod_{h<i} b_{h+p}^{\mathbbm{l}_{h}(\gamma_{1})+\mathbbm{l}_{h}(\gamma_{2})}\right) \left( b_{i+p}^{\mathbbm{l}_{i}(\gamma_{1})}-\E\left[ b_{i+p}^{\mathbbm{l}_{i}(\gamma_{1})}\right] \right)\left( b_{i+p}^{\mathbbm{l}_{i}(\gamma_{2})}-\E\left[ b_{i+p}^{\mathbbm{l}_{i}(\gamma_{2})}\right] \right) \left(\prod_{i<h} \E\left[b_{h+p}^{\mathbbm{l}_{h}(\gamma_{2})}\right]\E\left[ b_{h+p}^{\mathbbm{l}_{h}(\gamma_{1})}\right] \right) \\
& +\sum_{i_{2}<i_{1}<0}\left(\prod_{h<i_{2}} b_{h+p}^{\mathbbm{l}_{h}(\gamma_{2})+\mathbbm{l}_{h}(\gamma_{1})}\right)  \left( \left( b_{i_{2}+p}^{\mathbbm{l}_{i_{2}}(\gamma_{2})}-\E\left[ b_{i_{2}+p}^{\mathbbm{l}_{i_{2}}(\gamma_{2})}\right] \right)b_{i_{2}+p}^{\mathbbm{l}_{i_{2}}(\gamma_{1})}\right)\\* 
&\qquad \times\left(\prod_{i_{2}<h< i_{1}} \E\left[b_{h+p}^{\mathbbm{l}_{h}(\gamma_{2})}\right]b_{h+p}^{\mathbbm{l}_{h}(\gamma_{1})} \right) \left( b_{i_{1}+p}^{\mathbbm{l}_{i_{1}}(\gamma_{1})}-\E\left[ b_{i_{1}+p}^{\mathbbm{l}_{i_{1}}(\gamma_{1})} \right] \right) \left(\prod_{i_{1}<h} \E\left[b_{h+p}^{\mathbbm{l}_{h}(\gamma_{1})}\right]\E\left[ b_{h+p}^{\mathbbm{l}_{h}(\gamma_{2})}\right] \right)
\end{align*}
}
After taking expectation in this formula, from the independence of the entries, one arrives at
\begin{align*}
\cov(a_{\gamma_{1}+p},a_{\gamma_{2}+p})= 
\sum_{i<0}\left(\prod_{h<i}\E\left[ b_{h+p}^{\mathbbm{l}_{h}(\gamma_{1})+\mathbbm{l}_{h}(\gamma_{2})}\right]\right)
\cov\left( b_{i+p}^{\mathbbm{l}_{i}(\gamma_{1})},b_{i+p}^{\mathbbm{l}_{i}(\gamma_{2})}\right) \left(\prod_{i<h} \E\left[b_{h+p}^{\mathbbm{l}_{h}(\gamma_{2})}\right]\E\left[ b_{h+p}^{\mathbbm{l}_{h}(\gamma_{1})}\right] \right)
\end{align*}
from which, according to \eqref{ee1'}, it follows that 
\begin{equation}\label{e15}
\lim_{p\to\infty}\frac{1}{p^{\alpha k_{1}+\alpha k_{2}-2\epsilon}}\cov(a_{\gamma_{1}+p},a_{\gamma_{2}+p})=\sum_{i<0}\left(\prod_{h<i}m_{\mathbbm{l}_{h}(\gamma_{1})+\mathbbm{l}_{h}(\gamma_{2})}\right)C(\mathbbm{l}_{i}(\gamma_{1}),\mathbbm{l}_{i}(\gamma_{2}))\left(\prod_{i<h} m_{\mathbbm{l}_{h}(\gamma_{1})}m_{\mathbbm{l}_{h}(\gamma_{2})}\right).
\end{equation}
Let's point out that in the case $\epsilon=0$, one has $C(k,l)=m_{k+l}-m_{k}m_{l}$ and in this case the formula simplifies to the one given in \eqref{D}. 

Next we deal with the case in which there are flat steps in $\gamma_{1}$ and/or $\gamma_{2}$.  Thus we need only consider the case $\epsilon=\alpha$.  In the first place if only one of them has a flat step then we may assume that $\gamma_{1}$ has one flat step and $\gamma_{2}$ does not.  Then we write $a_{\gamma_{1}+p}=d_{g+p}\prod_{h<0}b_{h+p}^{\mathbbm{l}_{h}(\gamma_{1})}$ where $g$ is the level of the flat step.  Hence $a_{\gamma_{2}+p}=\prod_{h<0}b_{h+p}^{\mathbbm{l}_{h}(\gamma_{2})}$ and then because of the independence of the entries, 
\[
\cov\left(a_{\gamma_{1}+p},a_{\gamma_{2}+p}\right)=
\E[d_{g+p}]\cov(a_{\gamma_{1}'+p},a_{\gamma_{2}+p}),
\]
where $\gamma_{1}'$ is the path obtained from $\gamma_{1}$ by removing the flat step and gluing together the remaining parts.  As a consequence of the above, the path $(\gamma_{1}',\gamma_{2})$ does not have a flat step and for the  $\Gamma=(\gamma_{1}',\gamma_{2})$ we can use \eqref{e15}.  Taking the limit over $p\to\infty$ and keeping in mind that now the path $\gamma_{1}'$ has length $k_{1}-1$, using the previous step we get that 
\[
\lim_{p\to\infty}\frac{1}{p^{\alpha( k_{1}+ k_{2}-2)}}\cov\left(a_{\gamma_{1}+p},a_{\gamma_{2}+p}\right)=0
\]
The third situation here is the one in which both $\gamma_{1}$ and $\gamma_{2}$ have a flat step.  In this case we write $a_{\gamma_{1}+p}=d_{g_{1}+p}\prod_{h<0}b_{h+p}^{\mathbbm{l}_{h}(\gamma_{1})}$ and $a_{\gamma_{1}+p}=d_{g_{2}+p}\prod_{h<0}b_{h+p}^{\mathbbm{l}_{h}(\gamma_{2})}$ where $g_{1}$ and $g_{2}$ are the levels of the flat steps in  $\gamma_{1}$ and $\gamma_{2}$.  Denote by $\gamma_{1}'$ and $\gamma_{2}'$ the paths obtained by removing the flat steps.  In this case $\Gamma=(\gamma_{1}',\gamma_{2}')$ has no flat steps and the length of $\gamma_{1}'$ is $k_{1}-1$, while the length of $\gamma_{2}'$ is $k_{2}-1$.  Now a simple calculation gives
\[
\cov(a_{\gamma_{1}+p},a_{\gamma_{2}+p})=\E[d_{g_{1}+p}]\E[d_{g_{2}+p}]\cov(a_{\gamma_{1}'+p},a_{\gamma_{2}'+p})+\cov(d_{g_{1}+p},d_{g_{2}+p})\E[a_{\gamma_{1}'+p}]\E[a_{\gamma_{2}'+p}].
\]
Therefore, using \eqref{e15}, \eqref{dlim} and \eqref{e2}, and noting that $\mathbbm{l}_{i}(\gamma_{1}')=\mathbbm{l}_{i}(\gamma_{1})$ and similarly $\mathbbm{l}_{i}(\gamma_{2}')=\mathbbm{l}_{i}(\gamma_{2})$ we get (cf. Remark \ref{r1}),
\[
\lim_{p\to\infty}\frac{1}{ p^{\alpha (k_{1}+ k_{2}-2)}}\cov(a_{\gamma_{1}+p},a_{\gamma_{2}+p})=
\begin{cases}
\sigma_{d}^{2}\prod_{h<0}m_{\mathbbm{l}_{h}(\gamma_{1})}m_{\mathbbm{l}_{h}(\gamma_{2})}=\sigma_{d}^{2}m_{2}^{k+l-2} &\text{if}\: g_{1}=g_{2} \\
0 &\text{otherwise}.
\end{cases}
\]
Return with the results of Cases 1, 2 and 3 to \eqref{limU}.  In computing the limit of \eqref{f:1}, realize that we need to worry about only the case $j$ even and partitions $\Delta$ of pairs. Then $\Gamma(\mathbf{k}:\Delta)$
 is one-to-one with $\times_{\{i,j\}\in\Delta}\Gamma(k_{i},k_{j})$.  Returning to \eqref{P}, a moment of thinking gives that 
 \[
 \E\left[ (S_{n}(P))^{j} \right]=|\text{pairs of}\:\{1,2,\dots,j\}|\left( \sum_{k,l\ge1}w_{k}w_{l}D(k,l) \right)^{j}.
 \]
 Since $|\text{pairs of}\:\{1,2,\dots,j\}|=\frac{j!}{2^{j/2}(j/2)!}$, which are the even moments of the normal $N(0,1)$, the rest follows. \qedhere
\end{proof}

\begin{corollary}
Assume that for $0< \epsilon\le \alpha$ $b_{n}/n^{\alpha}=1+Z_{n}/n^{\epsilon}$, where $\lim_{n\to\infty}Z_{n}=Z$ is a random variable with finite moments, mean $0$ and variance $\sigma^{2}_{Z}$.  Then $C(k,l)=kl\sigma_{Z}^{2}$ and 
\[
D(k,l)= 
\begin{cases}
\frac{kl}{\alpha(k+l)}{k\choose k/2}{l\choose l/2} &\text{if}\:k,l\:\text{even}\\
\frac{kl}{\alpha(k+l)}{k-1\choose (k-1)/2}{l-1\choose (l-1)/2} &\text{if}\:\epsilon=\alpha\:\text{and}\:k,l\:\text{odd}\\
0 & otherwise.
\end{cases}
\]
In particular CLT holds for the model \eqref{e0}.  
\end{corollary}

\begin{proof}
To compute $C(k,l)$, just notice that (in moments)
\[
(b_{n}/n^{\alpha})^{k}\sim 1+kZ_{n}/n^{\epsilon}+O(1/n^{2\epsilon}),
\]
from which the formula of $C(k,l)$.

For the rest, there is only one thing we need to do, namely compute 
\[
\sum_{(\gamma_{1},\gamma_{2})\in\Gamma(k,l)}\sum_{i<0}\mathbbm{l}_{i}(\gamma_{1})\mathbbm{l}_{i}(\gamma_{2}).
\]
To carry this out, we fist realize the pairs $(\gamma_{1},\gamma_{2})\in\Gamma(k,l)$ by fixing $\gamma_{1}$ in $\Gamma_{k}$ and then ``sliding'' up and down another fixed $\zeta\in\Gamma_{l}$, to justify that  
\[
\sum_{i\in\Z}\sum_{\gamma_{1}\in\Gamma_{k}}\mathbbm{l}_{i}(\gamma_{1})\sum_{\zeta\in\Gamma_{l}}\sum_{q\in\Z}\mathbbm{l}_{i}(\zeta+q).
\]
For $i$ and $\zeta\in\Gamma_{l}$,  $\sum_{q\in\Z}\mathbbm{l}_{i}(\zeta+q)=l$ and since  there are ${l\choose l/2}$ paths in $\Gamma_{l}$, the rest follows. 
\end{proof}

\section{A Flavor of Free Probability Theory}\label{3}

Given independent tridiagonal matrices $A_{1,n}, A_{2,n},\dots,A_{l,n}$, one can ask about the joint distribution in the moment sense.  More precisely, is it true that (here $X_{u,n}=\frac{1}{n^{\alpha_{u}}}A_{u,n}$)
\[
\lim_{n\to\infty}\tr_{n}(X_{i_{1},n}X_{i_{2},n}\dots X_{i_{k},n})\:\text{exists for any}\:i_{1},i_{2},\dots,i_{m}?
\]
The answer is yes, but before we do that we need to give a definition.  We think about the set $\mathcal{C}=\{1,2,\dots, r \}$ as a set of colors.  Then, for a string of colors $\mathbf{c}=( i_{1},i_{2},\dots,i_{k} )$,  from the set $C$ we define 
\[
\Gamma_{k}^{\mathbf{c}}=\{ \gamma\in\Gamma_{k}:\text{each edge}\:j_{u},j_{u+1}\:\text{is colored with}\: i_{u} \}.
\]  
For a given color $u\in\mathcal{C}$, and a colored path $\gamma\in\Gamma_{k}^{\mathbf{c}}$, we define $\mathbbm{l}_{i}^{u}(\gamma)$ as the number of crossings of the line $i+1/2$ with the steps of $\gamma$ colored with $u$.  Similarly $\mathbbm{f}_{i}^{u}(\gamma)$ is the number of flat steps at level $i$ colored with color $u$.  

\begin{theorem}
Assume that for each $u\in\mathcal{C}$, the entries of the matrix $A_{u,n}$ satisfy, for some $\alpha_{u}>0$,
\[
\lim_{n\to\infty} \E \left[ \left(b_{u,n}/n^{\alpha_{u}}\right)^{k} \right] = m_{u,k}\quad\text{for any}\quad k\ge0
\]
with $m_{u,0}=1$ and 
\[
\sup_{n\ge1}\E\left[ \left|d_{u,n}\right|^{k} \right] < \infty\quad\text{for any}\quad k\ge0.
\]
Under these assumptions, if all the entries of the matrices are independent of one another and $\mathbf{c}=(i_{1},i_{2},\dots,i_{k})$, then
\[
\lim_{n\to\infty}\tr_{n}(X_{i_{1},n}X_{i_{2},n}\dots X_{i_{k},n})=\frac{1}{\alpha_{1}+\alpha_{2}+\dots+\alpha_{k}+1}\sum_{\gamma\in\Gamma_{k}^{\mathbf{c}}}\prod_{u\in\mathcal{C}}\prod_{i\in\Z}m_{u,\mathbbm{l}_{i}^{u}(\gamma)}
\]
where the limit is in expectation and also almost surely.  
\end{theorem}

One possible interpretation of this in term of noncommutative probability theory is the following.  Assume that $(\mathcal{X},\phi)$, $(\mathcal{Y},\psi)$ are noncommutative probability spaces, i.e. $\mathcal{X}$, $\mathcal{Y}$ are unital algebras over the complex numbers and $\phi:\mathcal{X}\to\C$,  $\psi:\mathcal{Y}\to\C$ are two linear functionals with $\phi(1)=\psi(1)=1$.  Assume $a_{1}',a_{2}',\dots a_{l}'$ are noncommutative random variables on $\mathcal{Y}$ such that $\psi((a_{u}')^{k})=m_{u,k}$.  Then the joint distribution of $a_{1},a_{2},\dots a_{l}$ is described by
\[\tag{*}
\bar{\phi}(a_{i_{1}}a_{i_{2}}\dots a_{i_{k}})=
\begin{cases} 
\frac{1}{\alpha_{1}+\alpha_{2}+\dots+\alpha_{k}+1}\sum_{\gamma\in\Gamma_{k}^{\mathbf{c}}}\prod_{u\in\mathcal{C}}\prod_{i\in\Z}\psi((a_{u}')^{\mathbbm{l}_{i}^{u}(\gamma)}) & k\:\text{even} \\
0 & k\:\text{odd}.
\end{cases}
\]
Note here that this dependence involves  $m_{u,k}$ for $k$ odd as well, as opposed to the defining relationship from \eqref{e2} which involved only $m_{k}$ for $k$ even.  

Take for example the case of just two such random matrices.  Rescale things out to have a nicer appearance to $\phi=(\alpha_{1}+\alpha_{2}+1)\bar{\phi}$. Let's take two random variables $a,b$.  Then, $\phi$ of a product of odd length in $a,b$ is $0$, while for products of even length we have (ignore here the presence of $\alpha_{1}$ and $\alpha_{2}$ or rescale the functional $\phi$),
\begin{align*}
\phi(a^{2})&=2\psi((a')^{2})\\
\phi(b^{2})&=2\psi((b')^{2})\\
\phi(ab)&=2\psi(a')\psi(b') \\
\phi(a^{4})&=2\psi((a')^{2})+4\psi((a')^{2})^{2}\\
\phi(b^{4})&=2\psi((b')^{2})+4\psi((b')^{2})^{2} \\ 
\phi(abab)&=2\psi((a')^{2})\psi((b')^{2})+4\psi(a')^{2}\psi(b')^{2} \\ 
\phi(a^{2}b^{2})&=4\psi((a')^{2})\psi((b')^{2})+2\psi(a')^{2}\psi(b')^{2} \\
\phi(a^{3}b)&= 2\psi((a')^{3})\psi(b')+4\psi((a')^{2})\psi(a')\psi(b').
\end{align*}
From this it's quite clear that, with respect to $\phi$, the moments of $a$ and $b$ alone do not determine their joint moments.  However, imposing the condition that $\psi(a^{2k+1})=\psi(b^{2k+1})=0$ for any $k\ge0$, one can do this.  For example in this case we have 
\begin{align*}
\phi(ab)&=0 \\ 
\phi(abab)&=\frac{1}{2}\phi(a^{2})\phi(b^{2})\\
\phi(a^{2}b^{2})&=\phi(a^{2})\phi(b^{2}) \\
\phi(a^{3}b)&=0. \\
\end{align*}

Another view at these things is the following.  Assume that $a'$ and $b'$ are independent random variables and the functional $\psi$ is just the expectation.  Then we consider sequences of iid random variables $\{ X_{i}\}_{i\in\Z}$ and $\{ Y_{i}\}_{i\in\Z}$ whose distributions are given by the distributions of $a'$ and $b'$.  Consider then the operators $\mathcal{A}$ and $\mathcal{B}$ given in \eqref{op}.  Now if $\mathcal{X}$ is the algebra of infinite dimensional random matrices like $\mathcal{A}$ and  the functional $\phi$ on the algebra generated by $\mathcal{A}$ and $\mathcal{B}$ is given by $\E[P(\mathcal{A},\mathcal{B})_{0,0}]$, for any noncommutative polynomial $P$ in two variables, then the joint distribution of $\mathcal{A}$ and $\mathcal{B}$ is given by (*).  

Returning to the general situation from (*), we want to point out that in the case that the variables $a_{u}'$ are symmetric, then the noncommutative joint moments of $a_{1},a_{2},\dots,a_{l}$ are given in terms of the individual moments of $a_{1},a_{2},\dots,a_{l}$.   This follows from the fact that all the $m_{k}$'s  involved in the joint moments have $k$ even and according to Remark~\ref{inv} these can be expressed back in terms of the moments of the variables $a_{i}$.  

We can call these variable ``independent'' in a certain way and interpret this fact via the relationship between the joint moments and the individual moments of each variable.  This can be seen by introducing some kind of cumulant and express this ``independence'' as a property of the cumulant.  In the classical or free cases of independence,t this corresponds to the simple fact that the joint cumulants are the sum of the cumulants of the individual variables.  

Finally if all the moments $m_{k}=1$, then the matrices $(X_{1,n},X_{2,n},\dots,X_{r,n})$ converges in distribution to $(S,S,\dots,S)$ where $S$ is a semicircular random variable, something not very interesting though but due to the fact that the coloring does not play any role here.  However if the moments $m_{k}$ are not constant equal to $1$, then the coloring does play an essential role.  

There is also a fluctuation result in this context as follows.  
\begin{theorem}\label{mdimf}
Assume that in addition to the properties in the above theorem we have that for some $0\le\epsilon_{u}\le\alpha_{u}$,
\[
\lim_{n\to\infty}n^{2\epsilon_{u}}\cov((b_{u,n}/n^{\alpha_{u}})^{k},(b_{u,n}/n^{\alpha_{u}})^{l})=C_{u}(k,l),
\]
\[
\E\left[\left|n^{\epsilon}((b_{u,n}/n^{\alpha_{u}})^{k}-\E[(b_{u,n}/n^{\alpha_{u}})^{k}])\right|^{m}\right]<\infty,\:\forall\: k,m\ge1,
\]
and 
\[
\lim_{n\to\infty}\var(d_{u,n})=\sigma_{u}^{2}.
\]
Now, take $\epsilon=\min_{u=1,\dots,r}(\epsilon_{u})$ and set 
\[
S_{n}(i_{1},i_{2},\dots,i_{k})=n^{\epsilon+1/2}\left(\tr_{n}(X_{i_{1},n}X_{i_{2},n}\dots X_{i_{k},n})-\E[\tr_{n}(X_{i_{1},n}X_{i_{2},n}\dots X_{i_{k},n})]\right).
\]
Then the family $\{ S_{n}(i_{1},i_{2},\dots,i_{k}):1\le i_{1},i_{2},\dots,i_{k}\le r \}$ converges to a Gaussian family.  
\end{theorem}

\section{Remarks and Extensions} \label{4}

\subsection{Still Tridiagonal Models}
There are various ways of extending Theorem~\ref{T1} and Theorem~\ref{fluct}. We refrained to give it in full because the proofs would have been overloaded with unnecessary notations and minor differences.  

In both theorems mentioned here we can replace the sequence $n^{\alpha}$ by any sequence $\alpha_{n}$ which satisfies the growth rate condition $\lim_{n\to\infty}n(\frac{\alpha_{n}}{\alpha_{n-1}}-1)=\alpha$.  

The second extension comes from allowing growth in the diagonal part.  Namely if we replace the condition \eqref{e1'} by the condition 
\[
\lim_{n\to\infty}\E[(d_{n}/n^{\beta})^{k}] = m_{k}'
\]
then, if $\beta<\alpha$, the same conclusion holds in Theorem~\ref{T1} and the same conclusion under Theorem~\ref{fluct} with the condition \eqref{dlim} replaced by $\lim_{n\to\infty}\var(d_{n}/n^{\beta})=\sigma_{d}$. However if $\beta=\alpha$, then the conclusions still hold, nonetheless the formulae of $L_{k}$ become
\[
L_{k}=\frac{1}{\alpha k +1}\sum_{\gamma\in\Gamma_{k}\cup\Gamma_{k}^{-}}\prod_{i\le0}m_{\mathbbm{l}_{i}(\gamma)}m_{\mathbbm{f}_{i}(\gamma)}',
\]
while in Theorem~\ref{fluct},  condition \eqref{dlim} has to be replaced by 
\[
\lim_{n\to\infty}n^{2\epsilon}\cov((d_{n}/n^{\alpha})^{k},(d_{n}/n^{\alpha})^{l})=C'(k,l)
\] 
 the only difference here is that the convariance matrix $D_{\epsilon}(k,l)$ now depends also on $C'(k,l)$ and $m_{k}'$.    

If $\beta>\alpha$, the scaling of the matrix $X_{n}$ has to be changed to 
$X_{n}=\frac{1}{n^{\beta}}A_{n}$.  The conclusions of both theorems hold with the appropriate changes since now the dominating terms are the ones on the diagonal.  For example the \eqref{e2}, becomes
\[
L_{k}=\frac{1}{\beta k+1}m_{k}'. 
\] 
We leave to the reader to see how the changes in Theorem~\ref{fluct} have to be done.  

Another extension is obtained by dropping the independence of the entries.  We can replace that in Theorem~\ref{T1} by a more relaxed version.  
\begin{remark}\label{nind}
Assume that for any $\gamma\in\Gamma_{k}\cup\Gamma_{k}^{-}$, there is a number $m_{\gamma}$ so that for a certain $\alpha>0$,
\[
\lim_{n\to\infty}\frac{1}{n^{k\alpha}}\E[a_{\gamma+n}]=m_{\gamma}.
\]
Then if $X_{n}=\frac{1}{n^{\alpha}}A_{n}$, 
\[
\lim_{n\to\infty}\tr_{n}(X_{n}^{k})=L_{k}
\]
where $L_{k}$ is computed by 
\[
L_{k}=\frac{1}{\alpha k+1}\sum_{\gamma\in\Gamma_{k}\cup\Gamma_{k}^{-}}m_{\gamma}.
\]
In particular one can apply this to the cases when the matrix $A$ is obtained from another tridiagonal matrix $B_{n}$ which has independent entries by replacing each entry with a function of the other nearby entries in a finite range.  For example one can replace the nonzero entries in $B_{n}$ by the average of the neighbors nearby it in a finite range.  Another example is the Laguerre $\beta$ models discussed in \cite{DE} and \cite{DE2}, or more general the models in which each entry in $B$ is replaced by a polynomial in the variables lying in finite neighborhood of the entry.  
\end{remark}

\subsection{Band Diagonal Models}

We can extend the results so far to a more general setting by allowing not only one subdiagonal but more than one.  In this case we take symmetric matrices of the form $A_{n}=\{ a_{i,j} \}_{i,j=1}^{n}$ so that $a_{i,j}=0$ for $|i-j|>w$, where $w$ is the width of the band and all entries are independent.  Denote $b_{v,i}=a_{i,v+i}$.  

In this case we can consider the problem of convergence of the eigenvalues and of the fluctuations. 
Before we give this extension, let us define the needed objects.  

Set 
\[
\Gamma_{k,w}=\{\gamma=(i_{1},i_{2},\dots,i_{k+1})\in Z^{k+1}:\:i_{1}=i_{k+1},\: |i_{u}-i_{u+1}|\le w,\max(\gamma)=0 \}.
\]
Then we define for any path $\lambda$, $\mathbbm{l}_{[i,j]}(\lambda)$ to be the number of steps $i_{u},i_{u+1}$ so that $\{i_{u},i_{u+1}\}=\{i,j\}$.  In particular, for the notations we already used we get $\mathbbm{l}_{i}(\lambda)=\mathbbm{l}_{[i,i+1]}(\lambda)$ and $\mathbbm{f}_{i}(\lambda)=\mathbbm{l}_{[i,i]}(\lambda)$.

Notice here the equivalent of the formula \eqref{ap} as
\[
a_{\lambda}=\prod_{i\in\Z,0\le v\le w}b_{v,i}^{\mathbbm{l}_{[i,i+v]}(\lambda)}.
\]

\begin{theorem}
Assume that for each $0\le v \le w$, and given $\alpha_{v}\ge0$, there is $m_{v,k}$ so that 
\[
\begin{cases}
\lim_{n\to\infty}\E[(a_{i,j}/n^{\alpha_{v}})^{k}]=m_{v,k}, &\:\text{if}\: |i-j|=v, \alpha_{v}=\alpha \\
\sup_{n}\E[(a_{i,j}/n^{\alpha_{v}})^{k}]<\infty, &\:\text{if}\: |i-j|=v, \alpha_{v}<\alpha \\
\end{cases}
\]
with $0<\alpha=\max(\alpha_{v}:1\le v \le w)$ and the convention that $m_{v,0}=1$ for any $v$.  Then, for $X_{n}=\frac{1}{n^{\alpha}}A_{n}$, one has that 
\[
\lim_{n\to\infty}\tr_{n}(X_{n}^{k})=L_{k}
\]
both in average and almost surely.  Moreover,
\[
L_{k}=\frac{1}{\alpha k+1}\sum_{\gamma\in\Gamma_{k,w}} \prod_{i\le 0,0\le v\le w}
\overline{m}_{\mathbbm{l}_{[i,i+v]}(\gamma)},
\]
where 
\[
\overline{m}_{k}=
\begin{cases}
m_{k} & \text{if}\: \alpha_{v}=\alpha\\
0 &\text{if}\:\alpha_{v}<\alpha.
\end{cases}
\]
\end{theorem}

This theorem says that in fact those subdiagonals not scaled by the maximum power $n^{\alpha}$, do not contribute to the limit $L_{k}$.  

Let us point out that one can extend this to a statement in which the independence condition is dropped and one gets a  version of Remark~\ref{nind}.

Similar versions of the first part of Proposition~\ref{ex} can be proved in this context too.  Namely, if each of the moments $(m_{v,k})_{k=1}^{\infty}$ come from the moments of a compactly supported measure, then the moments $L_{k}$ also come from a compactly supported measure.  In addition, if there are some numbers $m_{v}$ so that  $m_{v,k}=m_{v}$ for all $k\ge1$, then the corresponding measure with the moments given by $L_{k}$ is the distribution of $\sum_{v=1}^{w}m_{v}(z^{v}+z^{-v})$ under the Haar measure of the circle $S^{1}$. 

There is also a version of Theorem~\ref{fluct}.  
\begin{theorem}
In addition to the conditions given in the above Theorem, assume that for each $0\le v\le w$, there is $0\le \epsilon_{v}\le \alpha_{v}$ so that for any $k,l\ge1$
\[
\lim_{n\to\infty}n^{2\epsilon_{v}}\cov((b_{v,n}/n^{\alpha_{v}})^{k},(b_{v,n}/n^{\alpha_{v}})^{k})=C_{u}(k,l)
\]
and 
\[
\E\left[\left|n^{\epsilon_{v}}((b_{v,n}/n^{\alpha_{v}})^{k}-\E[(b_{v,n}/n^{\alpha_{v}})^{k}])\right|^{m}\right]<\infty,\:\forall\: k,m\ge1,
\]
Let $\epsilon=\min(\epsilon_{v}:0\le v\le w)$ and define 
\[
S_{n}(k)=n^{\epsilon+1/2}(\tr_{n}(X_{n}^{k})-\E[\tr_{n}(X_{n}^{k})]).  
\]
Then the family $\{ S_{n}(k)\}_{k=1}^{\infty}$ converges to a family of Gaussian random variables.  
\end{theorem}

\bibliography{Tridiagonal}
\bibliographystyle{amsplain}

\nocite{*}

\end{document}